\documentclass[11pt]{article}

\usepackage{amsmath,amssymb}
\usepackage{amsthm}
\usepackage[nobysame]{amsrefs}
\usepackage{indentfirst}
\usepackage{mathrsfs}
\usepackage{graphicx}
\usepackage{hyperref}
\usepackage{xcolor}
\usepackage{fourier}
\usepackage[margin=1.5in]{geometry}
\usepackage[norelsize,boxed,linesnumbered,vlined,algo2e]{algorithm2e} 

\usepackage{mathtools}
\makeatletter
\DeclareRobustCommand*\cal{\@fontswitch\relax\mathcal}
\makeatother

\usepackage{booktabs}

\usepackage{latexsym,enumerate}
\usepackage{multirow}

\newcommand{\comment}[1]{}
\newcommand{\EEA}{\end{eqnarray}}
\newcommand{\BEA}{\begin{eqnarray}}
\newcommand{\peq}{p_{eq}}

\begin{document}

\title{Diffusion forecasting model with basis functions from QR decomposition}

\author{John Harlim$^{1,2}$\footnote{jharlim@psu.edu} and Haizhao Yang$^{3}$\footnote{matyh@nus.edu.sg}\\
  \vspace{0.1in}\\
  $^{1}$ Department of Mathematics, the Pennsylvania State University,\\
   University Park, PA 16802-6400, USA.\\
  $^{2}$ Department of Meteorology and Atmospheric Science,\\ the Pennsylvania State University, University Park, PA 16802-5013, USA.\\
  $^{3}$ Department of Mathematics, National University of Singapore,\\
   10 Lower Kent Ridge Road, Singapore
}

\maketitle

\begin{abstract}
The diffusion forecasting is a nonparametric approach that provably solves the Fokker-Planck PDE corresponding to It\^o diffusion without knowing the underlying equation. The key idea of this method is to approximate the solution of the Fokker-Planck equation with a discrete representation of the shift (Koopman) operator on a set of basis functions generated via the diffusion maps algorithm. While the choice of these basis functions is provably optimal under appropriate conditions, computing these basis functions is quite expensive since it requires the eigendecomposition of an $N\times N$ diffusion matrix, where $N$ denotes the data size and could be very large. For large-scale forecasting problems, only a few leading eigenvectors are computationally achievable. To overcome this computational bottleneck, a new set of basis functions constructed by orthonormalizing selected columns of the diffusion matrix and its leading eigenvectors is proposed. This computation can be carried out efficiently via the unpivoted Householder QR factorization. The efficiency and effectiveness of the proposed algorithm will be shown in both deterministically chaotic and stochastic dynamical systems; in the former case, the superiority of the proposed basis functions over purely eigenvectors is significant, while in the latter case forecasting accuracy is improved relative to using a purely small number of eigenvectors. Supporting arguments will be provided on three- and six-dimensional chaotic ODEs, a three-dimensional SDE that mimics turbulent systems, and also on the two spatial modes associated with the boreal winter Madden-Julian Oscillation obtained from applying the Nonlinear Laplacian Spectral Analysis on the measured Outgoing Longwave Radiation (OLR). 
\end{abstract}

\maketitle

\section{Introduction}
Probabilistic prediction of dynamical systems is an important subject in applied sciences with wide applications, e.g., numerical weather prediction \cite{ehrendorfer:94}, population dynamics \cite{hbda:09}, etc. The ultimate goal of the probabilistic prediction is to provide an estimate of the future state as well as the uncertainty of the estimate. Such information is usually encoded in the probability densities of the state variables at the future times. 

In the context of stochastic dynamical systems, the forecasting problem can be illustrated as follows. Let $x(t)\in\mathcal{M}\subset \mathbb{R}^n$ be the state variable that satisfies,
\begin{align}\label{SDE} dx = a(x)\, dt + b(x)\, dW_t, \quad x(0)=x_0\sim p_0(x),\end{align}
where $W_t$ is the standard Wiener process, $a(x)$ is the vector field, and $b(x)$ is the diffusion tensor, all of which are defined on the manifold $\mathcal{M}\subset\mathbb{R}^n$ such that the solutions of this SDE exist. Assume that the stochastic process in \eqref{SDE} can be characterized by density functions $p(x,t)$ that satisfies the Fokker-Planck equation,
\begin{align}
\frac{\partial p}{\partial t} = \nabla\cdot (-ap + \frac{1}{2}\nabla (bb^\top)p) := \mathcal{L}^*p , \quad p(x,0) = p_0(x).\label{FokkerPlanck}
\end{align}
Here the differential operators are defined with respect to the Riemannian metric\index{Riemannian metric} inherited by $\mathcal{M}$ from the ambient space $\mathbb{R}^n$. Assume that the dynamical system in \eqref{SDE} is ergodic such that $p(x,t)\to\peq(x)$ as $t\to\infty$, where $\peq$ denotes the equilibrium density that satisfies $\mathcal{L}^*\peq = 0$. This ergodicity assumption also implies that the manifold $\mathcal{M}$ is the attractor of the dynamical system in \eqref{SDE}. In this setting, the goal of the  probabilistic prediction is to solve \eqref{FokkerPlanck} and compute the statistical values,
\BEA
\mathbb{E}[A(X)](t)  = \int_\mathcal{M} A(x) p(x,t) dV(x),\nonumber
\EEA
for any functional $A(x) \in L^1(\mathcal{M},p)$ where $V(x)$ denotes the volume form inherited by $\mathcal{M}$ from $\mathbb{R}^n$. 

If $a(x)$ and $b(x)$ are known, classical approaches for solving this problem are to solve the PDE in \eqref{FokkerPlanck} or to approximate $p(x,t)\approx \sum_k \delta(x-x_k(t))$, where $x_k(t)$ is the solution of the SDE in \eqref{SDE} given initial conditions $x_k(0)$. This latter technique, which is more feasible for high-dimensional problems, is often coined as Monte-Carlo \cite{leith:74} or ensemble forecasting \cite{tk:97}. When coefficients $a(x)$ and $b(x)$ are not available, a classical nonparametric approach is to estimate these coefficients without assuming any specific functional form \cite{PRL1,PRL2}. Such an approach uses the Kramer-Moyal expansion to approximate these coefficients by averaging over $\mathbb{R}^n$, which is impractical when $n$ is large. {\color{black} One way to overcome this issue is to seek for the intrinsic low-dimensional coordinate, e.g., using \cite{cao:97}, prior to the nonparametric fitting.} A more elegant technique is the so-called {\it diffusion forecast}, which makes no attempt to estimate these coefficients \cite{bgh:15}. Instead, it approximates the semigroup operator, $e^{\tau\mathcal{L}}$, of the generator $\mathcal{L}$ that is adjoint to $\mathcal{L}^*$ in \eqref{FokkerPlanck} with a set of data-adapted orthogonal basis functions defined on the intrinsic manifold $\mathcal{M}$. The advantage of this approach is that it is numerically feasible even if the ambient space $\mathbb{R}^n$ is high-dimensional as long as the intrinsic manifold $\mathcal{M}$ is low-dimensional. Given data samples $\{x_i\}_{i=1}^N$ with sampling density $\peq$, this method represents the probability density function $p(x,t)$ as a linear superposition of basis functions coming from the eigendecomposition of an $N\times N$ matrix $T$, which is a discrete representation of the operator $e^{\epsilon\hat{\cal L}}$, where $\hat{\cal L}: = \peq^{-1}\mbox{div}(\peq\nabla\, )$, is constructed with the diffusion maps algorithm \cite{cl:06} with bandwidth parameter $\epsilon$. The representation of the probability densities with these basis functions was proven to be optimal under appropriate condition \cite{bgh:15}, and the resulting forecasting model provably solves the Fokker-Planck PDE in \eqref{FokkerPlanck} (see Chapter~6 of \cite{harlim:18}); these results have been numerically verified in various applications \cite{bgh:15,bh:16physd,bh:16jcp}. However, this approach requires a large number of basis functions for accurate solutions. Practically, it is often the case that only a few leading eigenvectors of $T$ can be numerically estimated accurately. The main goal of this article is to introduce a new set of basis functions that can be computed efficiently. This new basis functions are generated from the QR-decomposition of a matrix consisting of selected columns of $T$ and it's few leading eigenvectors. We will show numerical results comparing the diffusion forecasting method produces with this basis representation with those that use purely eigenvectors and QR-decomposition of purely columns of $T$. 

The remainder of this article is organized as follows. In Section~\ref{sec2}, we review the diffusion forecasting method. In Section~\ref{sec3}, we present the basic idea of the new basis functions. In Section~\ref{sec4}, we discuss two methods to specify initial conditions for the diffusion forecasts and present numerical results on both deterministic and stochastic dynamical systems. In Section~\ref{sec5}, we test our approach on the two spatial modes associated with the boreal winter Madden-Julian Oscillation. In Section~\ref{sec6}, we conclude the paper with a short discussion. 

\section{Diffusion forecast}\label{sec2}

Let $\{\varphi_k\}$ be the eigenfunctions of the weighted Laplacian operator $\hat{\cal L}: = \peq^{-1}\mbox{div}(\peq\nabla\, )$. Assuming that $\hat{\cal L}$ has a pure point spectrum, then these eigenfunctions form an orthonormal basis of $L^2(\mathcal{M},\peq)$. The key idea of the diffusion forecasting method \cite{bgh:15} is to represent the solutions of \eqref{FokkerPlanck}, $p(x,t)=e^{t\mathcal{L}^*}p_0(x)$, as follows,
\BEA
p(x,t) = \sum_{k} c_k(t) \varphi_k(x) \peq(x),\label{fouriersum}
\EEA 
such that,
\BEA
c_k(\tau) &=&    \langle e^{\tau\mathcal{L}^*}p_0, \varphi_k \rangle = \langle p_0, e^{\tau\mathcal{L}}\varphi_k \rangle =   \langle \sum_{j=0}^\infty c_j(0)\varphi_j\peq , e^{\tau\mathcal{L}}\varphi_k \rangle \nonumber \\ &=&\sum_{j=0}^\infty \langle \varphi_j, e^{\tau\mathcal{L}}\varphi_k \rangle_{\peq}  c_j(0).\label{DF:c}
\EEA
In the second equality above, we have used the fact that $e^{\tau\mathcal{L}}$ is adjoint with respect to $e^{\tau\mathcal{L}^*}$ and in the third equality above, we have used the expansion in \eqref{fouriersum} for the initial condition, $p_0(x)$. Computationally, we will apply a finite summation up to mode $M-1$. Obviously, this means that Eq.~\eqref{DF:c} is simply a matrix vector multiplication $\vec{c}(\tau) = A\vec{c}(0)$, where the $k$th component of $\vec{c}(\tau)\in\mathbb{R}^{M}$ is $c_k(\tau)$, and the components of matrix $A\in\mathbb{R}^{M\times M}$ are 
\BEA
A_{kj} = \langle \varphi_j, e^{\tau\mathcal{L}}\varphi_k \rangle_{\peq}.\label{DF:matrixA} 
\EEA
For $t=n\tau$, we can iterate the matrix $A$ to obtain $\vec{c}(t) = A^n\vec{c}(0)$. 

The second crucial idea in the diffusion forecasting is to approximate the matrix $A$ since the Fokker-Planck operator $\mathcal{L}^*$ is unknown. In order to estimate the components of $A$, we need to approximate $e^{\tau\mathcal{L}}\varphi_k$ in \eqref{DF:matrixA}. \index{Dynkin's formula} From the Dynkin's formula \cite{oksendal:03}, the solutions of the backward Kolmogorov equation with initial condition $f(x_i)\in\mathcal{C}^2(\mathcal{M})$ for the It\^o diffusion in \eqref{SDE} can be expressed as follows,
\BEA
e^{\tau\mathcal{L}}f(x_i) = \mathbb{E}_{x_i} [f(x_{i+1})],\label{DF:dynkin}
\EEA   
where $\mathbb{E}_{x_i}$ denotes the expectation conditional to the state $x_i$ and we define $\tau=t_{i+1}-t_i>0$. Define the shift operator $S_\tau$ (also known as the Koopman operator when $x$ is deterministic) on functions $f\in L^2(\mathcal{M},\peq)\cap C^2(\mathcal{M})$ as $S_\tau f(x_i)=f(x_{i+1})$ and substitute this equation to \eqref{DF:dynkin}, we have,
\BEA
e^{\tau\mathcal{L}}f(x_i) = \mathbb{E}_{x_i} [S_\tau f(x_{i})].\label{shiftapprox}
\EEA   
Equation \eqref{shiftapprox} suggests that $S_\tau$ is an unbiased estimator for $e^{\tau\mathcal{L}}$. Therefore, the components of matrix the $A$ can be approximated as follows,
\BEA
A_{kj} = \langle \varphi_j, e^{t\mathcal{L}}\varphi_k \rangle_{\peq} \approx \langle \varphi_j, S_\tau\varphi_k \rangle_{\peq} = \tilde{A}_{kj}.\label{DF:approxA1}\EEA
Numerically, we can estimate $\tilde{A}_{kj}$ with a Monte-Carlo integral,
\BEA
\tilde{A}_{kj} \approx \hat{A}_{kj} := \frac{1}{N-1}\sum_{i=1}^{N-1} \varphi_j(x_i) \varphi_k(x_{i+1}).\label{DF:approxA}
\EEA
With this approximation, notice that the solutions of the Fokker-Planck equation in \eqref{FokkerPlanck} are approximated with the representation in \eqref{fouriersum}, \eqref{DF:c}, \eqref{DF:matrixA}, and \eqref{DF:approxA}. 

Notice that this approximation is nonparametric since all we need are the basis functions $\varphi_k$ of $L^2(\mathcal{M},\peq)$ which are the eigenfunctions of the operator $\hat{\cal L}$ (and also eigenfunctions of $e^{\epsilon\hat{\cal L}}$). {\color{black} In numerical implementation, the forecasting accuracy can be improved by minimizing the approximation errors in the finite summation approximation to \eqref{fouriersum}-\eqref{DF:c}, in the unbiased estimation in \eqref{DF:approxA1}, and in the Monte-Carlo integral in \eqref{DF:approxA}. Though the eigenfunctions of the operator $\hat{\cal L}$ can be proved to be the optimal choice for the approximation in \eqref{DF:approxA1}, they might not be optimal for the approximation in \eqref{fouriersum}-\eqref{DF:c}, especially when the density function $p(x,t)$ contains sharp changes, in which case we expect that local basis functions could be better if they can significantly reduce the approximation error in \eqref{fouriersum}-\eqref{DF:c} while keeping the approximation error in \eqref{DF:approxA1} reasonably small. This observation motivates our work in the next section to explore other choices of orthogonal basis functions of the range space of the operator $e^{\epsilon\hat{\cal L}}$.}

\section{Empirical basis functions from QR decomposition}\label{sec3}
Before discussing the new basis functions that are obtained via the QR algorithm, we first provide a brief review of the diffusion maps construction in approximating $\mathcal{\hat{L}}$.
Given a fixed bandwidth Gaussian kernel
\BEA
K_\epsilon(x,y) =exp\Big(-\frac{\|x-y\|^2}{4\epsilon}\Big),\nonumber
\EEA
where $\epsilon$ denotes the bandwidth parameter, and a data $\{x_i\}_{i=1,\ldots,N}$, the diffusion maps algorithm constructs an $N\times N$ matrix $K_{ij}$ whose components are
\BEA
K_{ij} := \hat{K}_\epsilon(x_i,x_j) := \frac{K_\epsilon(x_i,x_j)}{q_\epsilon(x_i)^\alpha q_\epsilon(x_j)^\alpha},\nonumber
\EEA
where $q_\epsilon(x_i) = \frac{1}{N}\sum_{j=1}^N K_\epsilon(x_i,x_j)$. For $\alpha=1/2$, if we define an $N\times N$ diagonal matrix $D$ with components $D_{ii} := \hat{q}(x_i) =  \frac{1}{N}\sum_{j=1}^N \hat{K}_\epsilon(x_i,x_j)$, then the matrix 
\BEA
L = \frac{1}{\epsilon}(D^{-1}K - I_N) \nonumber
\EEA
converges to $\mathcal{\hat{L}}$ as $N\to\infty$ (see e.g., the Appendix of \cite{bh:16physd} or Chapter~6 of \cite{harlim:18} for the detailed derivations). In our numerical implementation, we will specify $\epsilon$ using the auto-tuning algorithm proposed in \cite{bh:16vb}. Here, the notation $I_N$ denotes the $N\times N$ identity matrix. 
Since $L$ is not symmetric, it is convenient to solve the eigenvalue problem of a conjugate symmetric matrix
\BEA
\hat{L} = D^{1/2} L D^{-1/2} = \frac{1}{\epsilon}(D^{-1/2}KD^{-1/2} - I_N)
\label{conjugation}
\EEA
such that $\hat{L}U=U\Lambda$, where the columns of $U$ are the eigenvectors of $\hat{L}$. Given these eigenvectors, one can verify that the eigenvectors of $L$ are the columns of $\Phi:= D^{-1/2}U$. Numerically, solving the eigenvalue problem for $\hat{L}$ directly may not be feasible especially when $\epsilon$ is small. Since $\hat{\cal L}$ and $e^{\epsilon\hat{\cal L}}$ share the same eigenfunctions, we usually solve a symmetric eigenvalue problem $\hat{T}U=US$, where, 
\BEA
\hat{T}= I_N + \epsilon \hat{L}.\label{T_hat}
\EEA 
Subsequently, the eigenvectors of 
\begin{equation}
\label{eqn:T}
T = I_N + \epsilon L
\end{equation}
are obtained via the conjugate relation $D^{-1/2}U=\Phi := [\vec\varphi_{1},\ldots,\vec\varphi_{N}]$. Here, the matrix $T$ is a discrete approximation of $e^{\epsilon\hat{\cal L}}$. So, the eigenvectors, $\vec{\varphi}_k$, of the matrix $T$ are a discrete approximation to the eigenfunctions, $\varphi_k(x)$, of the operator $e^{\epsilon\mathcal{\hat L}}$. More specifically, the $i$th component of $\vec{\varphi}_k$ is a discrete approximation to the function value $\varphi_k(x_i)$, which is the eigenfunction $\varphi_k(x)$ evaluated at the training data $x_i$. {\color{black} In general, the construction of $T$ requires no particular ordering of data. For the diffusion forecast application, however, it is more convenient to order the rows and columns of $T$ in accordance to the time series since the construction of \eqref{DF:approxA} requires the knowledge of temporal ordering.}

One of the practical issues with the representation of density $p$ using the basis functions in \eqref{fouriersum} is that when finite number of modes, $M$, is used, it may require $M$ to be large enough to avoid the Gibbs phenomenon, which causes inaccurate estimation of the normalization constant of $p$ and its statistics, especially when $p$ contains local events. In this case, it is standard that in approximation theory basis functions with local supports are used to represent the density $p$; this motivates the local basis functions in this paper. From the computation point of view, although many advanced algorithms have been proposed recently, computing a large portion of the eigenvectors of the matrix $T$ is still impractical when $N$ is large. Even if the eigendecomposition is performed in a distributed and parallel computer cluster, the scalability of modern eigensolvers is still limited due to the expensive communication cost. Hence, in practice, only a few leading eigenfunctions corresponding to the eigenvalues close to $1$ is available. Our goal is to propose a practical algorithm that computes a new set of basis functions for large scale diffusion forecasting, utilizing the available $M_E$ leading eigenfunctions with $M_E\ll N$. The space spanned by this new basis functions is a superset of the space spanned by these $M_E$ eigenfunctions.

From the matrix standpoint, the eigenvectors of $T$ form an orthogonal basis for the range of $T$, that is, $\vec{\varphi}_k^\top  \vec{\varphi}_\ell = N \delta_{k,\ell}$, where $\delta_{k,\ell}=1$ if $k=\ell$ and zero otherwise. Our key idea to exploit the fact that the range of $T$ has non-unique sets of orthogonal basis and to employ the $QR$-decomposition as a possible method to construct these basis functions. In particular, we would like to construct a subspace of the range of $T$ that is also a superset of the subspace spanned by the available leading eigenvectors of $T$ as follows. For the purpose of convenience, we will use the notation in MATLAB in Algorithm \ref{alg:qr}. Without loss of generality, we assume that the size of the data, $N$, and the number of selected columns from $T$, $M_Q$ are even in Algorithm \ref{alg:qr}.

\begin{algorithm2e}[H]
\label{alg:qr}
\caption{Mixed mode algorithm for constructing orthonormal basis. }
Input: A matrix $T= I_N + \epsilon L$ of size $N\times N$, and its first $M_E$ leading eigenvectors forming a matrix $\Phi$, and $M_Q$ as the number of columns selected from $T$.

Output: $M:=M_E+M_Q$ orthonormal basis vectors $\{ \vec\varphi_{j} \}_{j=1,\dots,M}$.

Let $T_{sub}:=T(:,\frac{N-M_Q}{2}:\frac{N+M_Q}{2}-1)$.

Let $B=[T_{sub},\Phi]$.

Compute the QR decomposition of $B$ via $[Q,R]=qr(B)$, where $qr(\cdot)$ is the truncated and unpivoted Householder QR decomposition returning $Q$ of size $N\times (M_E+M_Q)$ and $R$ of size $(M_E+M_Q)\times (M_E+M_Q)$.

The columns of $Q$ form $M_E+M_Q$ orthonormal basis vectors $\{ \vec\varphi_{j} \}_{j=1,\dots,M_E+M_Q}$.
\end{algorithm2e}
\vspace{10pt}

\begin{figure}
  \centering
  \includegraphics[width=.9\linewidth,height=12cm]{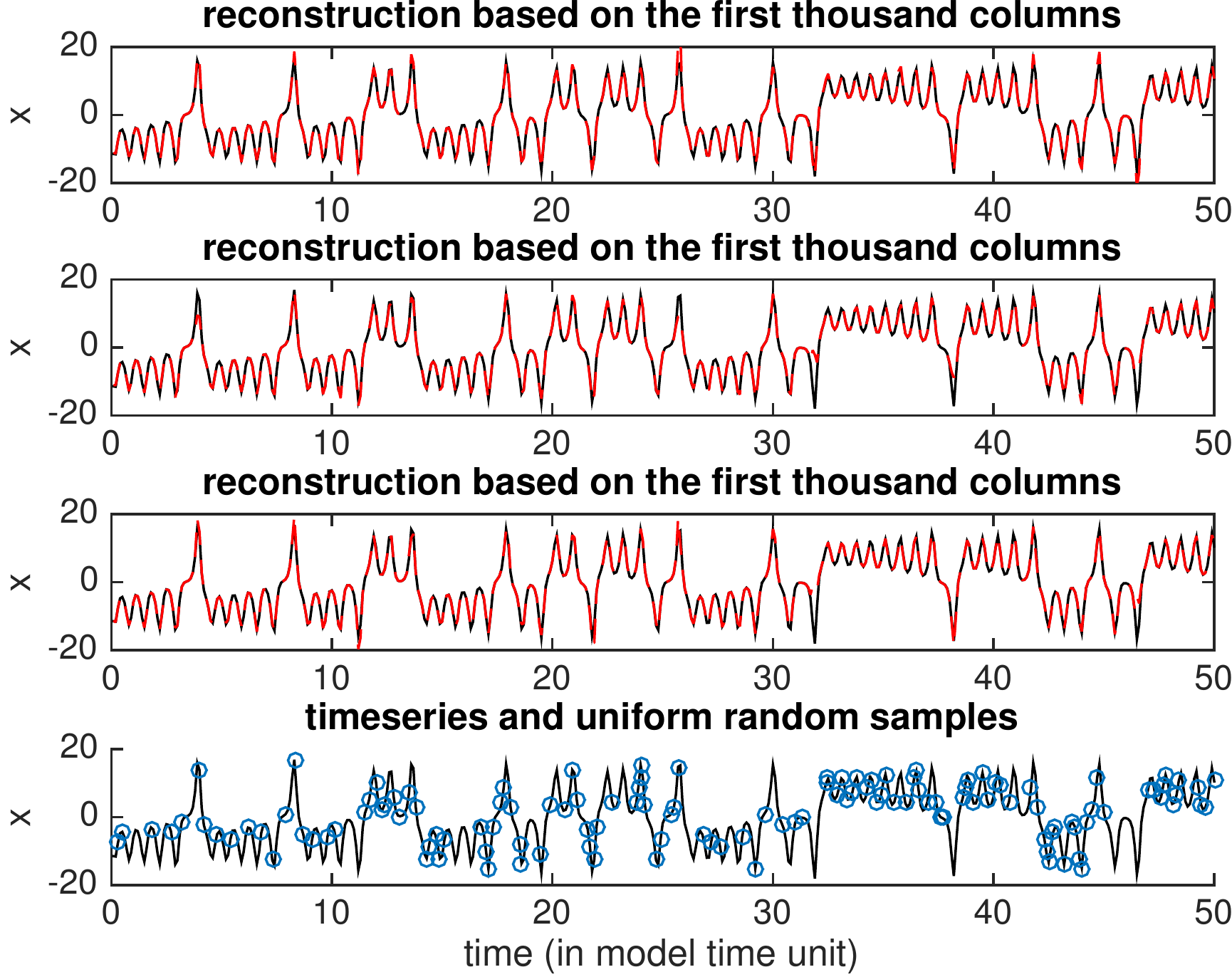}
  \caption{In each column, we show the true time series $x$ (black) on training data time interval $[0,50]$ in model time unit. The reconstructed time series (red) are based QR-decomposition of $T_{sub}$, chosen based on the first thousand columns of $T$ (first row), the second thousand columns of $T$ (second row), and 1000 uniformly distributed randomly selected columns of $T$ (third row). The randomly chosen columns of $T$ for the result above correspond to the data at time indices labeled with the blue circles (see fourth row).}
\label{L63recon} 
\end{figure}

{\color{black} We should point out that each column of $T$ contains the affinity between the training data set to the data corresponding to this column; a nonzero component implies that these pairs of data have stronger affinity and vice versa. Since the data corresponding to each column is a particular realization (or sample) of the invariant measure of the dynamics, this suggests that the columns of $T_{sub}$ can be chosen to be any arbitrary $M_Q$ sub-columns of $T$ so long as $M_Q$ is large enough to sample the invariant measure of the dynamics. To support this empirical argument, we show reconstructions of the training data using different sets of basis functions, obtained by $QR$-decomposition of various choices of $T_{sub}$. In particular, we consider a training data set of size $N=4000$ from the chaotic three-dimensional Lorenz-63 model (see Section~4.3 for the governing equations) such that the resulting diffusion matrix is $T\in\mathbb{R}^{4000\times4000}$. In Figure~\ref{L63recon}, we show the reconstruction of the first component, $x$, on time interval $[0, 50]$ using the basis functions obtained via Algorithm~1 with $M_E=0$ and $T_{sub}$ are chosen three different ways: the first thousand columns of $T$ (first row), the second thousand column of $T$ (second row), and 1000 uniformly distributed randomly selected columns of $T$ (third row). The randomly chosen columns of $T$ for the results above correspond to the data at time indices labeled in blue circles (see the fourth row of Figure~\ref{L63recon}).  From this numerical experiment, we found that the qualities of the reconstructions are relatively similar. Closer inspection indicates that one reconstruction can be better than the other on one regime but worse on the other regime. Intuitively, a better reconstruction should depends on whether there are enough data to represent the dynamical regimes. So, the choice of $T_{sub}$ as in Step 3 in Algorithm~\ref{alg:qr}, that is, the consecutive $M_Q$ middle sub-columns of $T$, is somewhat arbitrary. We prefer this choice with large enough $M_Q$ over the uniformly distributed randomly chosen columns to avoid undersampling part of the attractor with high value of sampling (or invariant) distribution. One possible way to avoid this issue is to choose the columns randomly in accordance to the sampling distribution of the data which we will not pursue in this article.}

In Algorithm \ref{alg:qr}, {\color{black} when $M_E$ is very small}, the most expensive part is the truncated and unpivoted Householder QR decomposition in Line $5$ with a computational complexity $O((M_E+M_Q)^2N)$. This QR decomposition can be implemented via Level-3 Basic Linear Algebra Subprograms
(BLAS) \cite{Dongarra:1990}, the operations of which can be implemented to achieve very high performance on modern processors \cite{Anderson:1999} and good scalability in distributed architectures \cite{Poulson:2013}. Comparing to computing $M_E+M_Q$ eigenvectors of $T$, Algorithm \ref{alg:qr} is highly efficient and admits the application to large scale diffusion forecasting problems.

In numerical linear algebra, eigendecomposition and rank-revealing QR decomposition are two standard tools to compress a matrix $T$. By selecting the leading eigenvectors or the important columns of $T$ according to the rank-revealing QR decomposition, one can obtain good approximations to the range of $T$. However in Algorithm \ref{alg:qr}, $M_Q$ consecutive columns of $T$ are selected for generating a set of orthonormal basis, instead of computing the rank-revealing QR factorization via column pivoting to choose $M_Q$ important columns of $T$. First, rank-revealing QR factorization with a large rank is more expensive and has a worse scalability in high-performance computing. Second, non-consecutive columns of $T$ coming from the rank-revealing QR decomposition correspond to time series that does not even fill up the attractor of the dynamical system and would result in poor forecasting accuracy.

\begin{figure}
  \centering
  \includegraphics[width=.49\linewidth,height=4cm]{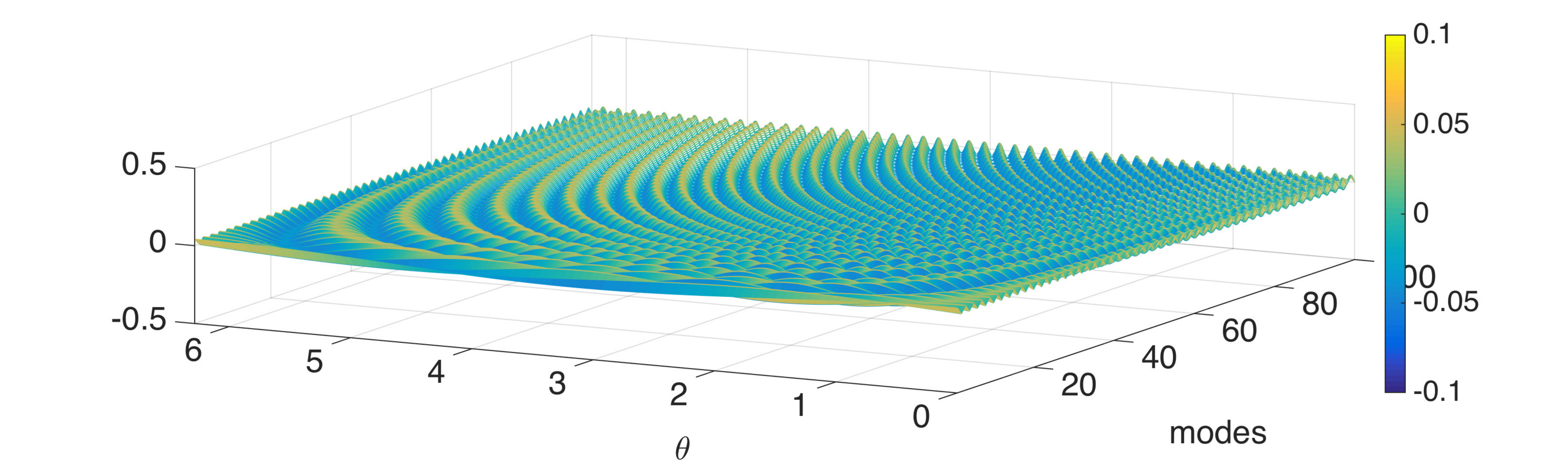}
  \includegraphics[width=.49\linewidth,height=4cm]{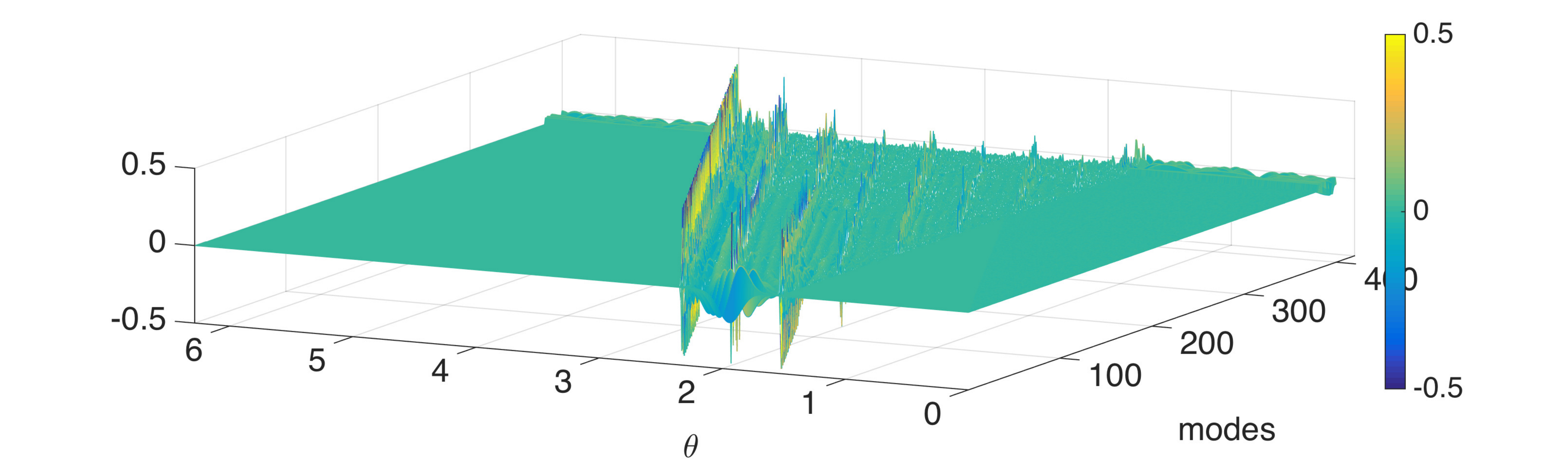}\\
  \includegraphics[width=.49\linewidth,height=4cm]{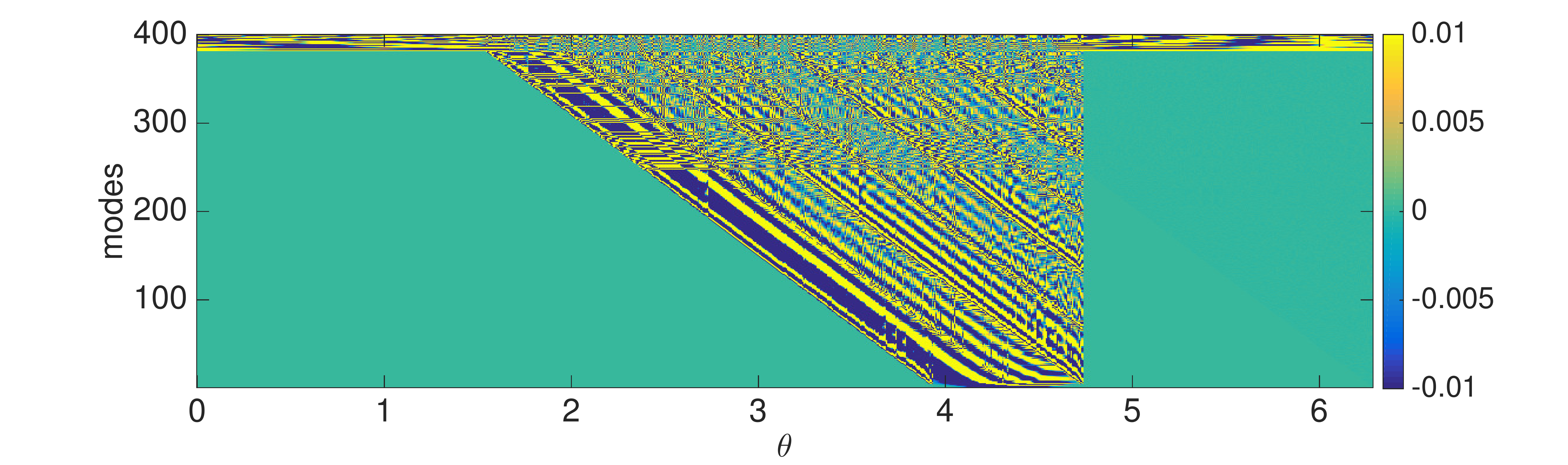} 
  \includegraphics[width=.49\linewidth,height=4cm]{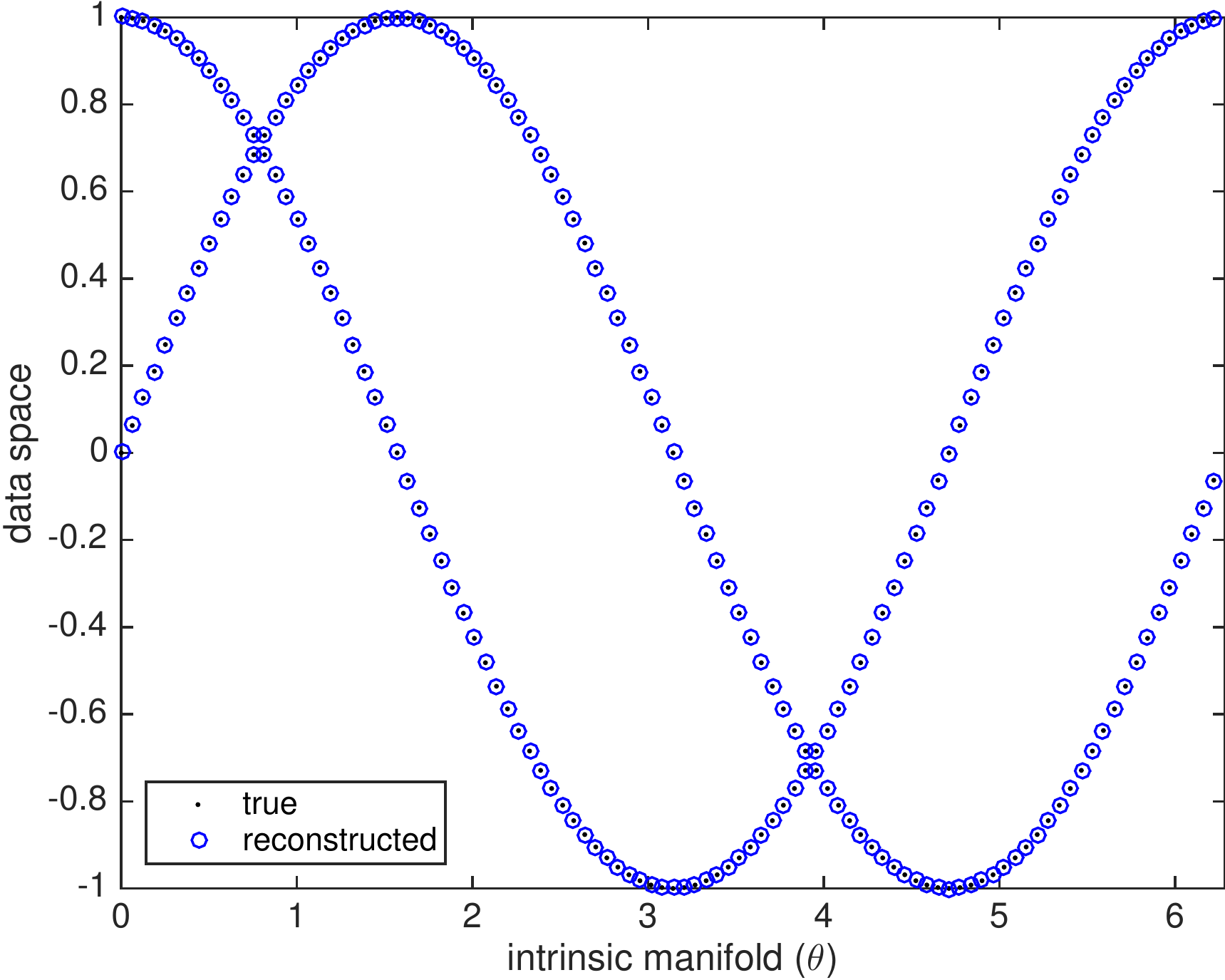}
  \caption{Top left: first $100$ leading eigenfunctions as a function of the intrinsic variable $\theta$ on periodic domain ${\cal M}=[0,2\pi)$. top right: $400$ basis functions by Algorithm \ref{alg:qr} with $M_E=20$ and $M_Q=380$. Bottom left: the same basis functions as in the middle panel. The color scale is set to $[-0.01,0.01]$ for the purpuse of better visualization. The sizes of the supports of the basis functions grow gradually. The basis function oscillates and behaves like a windowed Fourier mode. Bottom right: The true data $x_i = (\cos(\theta_i),\sin(\theta_i))$ and the reconstructed data using 400 basis functions by Algorithm \ref{alg:qr} with $M_E=20$ and $M_Q=380$. For clarity of presentation, we only plot every 10 data point from the total of equally spaced $N=1000$ training data points. The error $\|x - QQ'x \|_2 = 5\times 10^{-6}$.}
\label{fig:basis} 
\end{figure}

As we shall see later in the numerical examples, the forecasting accuracy using mixed modes modes in Algorithm \ref{alg:qr} is better than the one using even $M_E+M_Q$ eigenfunctions in deterministic dynamics. When the initial condition of Equation \eqref{FokkerPlanck} has local events, the localized basis from the QR factorization of the sparse matrix $T$, as oppose to the global basis from the eigendecomposition of $T$, would lead to better approximation accuracy when we numerically implement the Galerkin method in \eqref{fouriersum}. When $M_Q$ is sufficiently large, $M_E$ can even be set to $0$ since the QR decomposition of the middle portion of $T$ has generated basis with large enough supports over the invariant measure of the dynamical system (See Figure \ref{fig:basis} (bottom)). 


To illustrate the basis functions constructed by Algorithm \ref{alg:qr}, let us take the example when the intrinsic manifold $\mathcal{M}$ is a unit circle and the observation data are uniformly sampled on $\mathcal{M}$ with $N=1000$. In this case, the eigenfunctions of $T$ are simply the Fourier modes as visualized in Figure \ref{fig:basis} (top). Algorithm \ref{alg:qr} is applied with $M_E=20$ and $M_Q=380$ to construct $400$ basis functions. As shown in Figure \ref{fig:basis} (top right and bottom left), the supports of the basis functions from mixed modes gradually grow from a local regime to the whole domain; the basis function oscillates within its supports behaving like a windowed Fourier mode. Hence, applying the basis functions from Algorithm \ref{alg:qr} to the Galerkin method in Equation \eqref{fouriersum} is essentially performing a multiresolution analysis on the density $p$. When $p$ contains local events, the basis functions from Algorithm \ref{alg:qr} can capture these features better than the eigenfunctions of $T$ that are all global. Therefore, in the case of local events, the basis functions from Algorithm \ref{alg:qr} would lead to better forecasting accuracy. {\color{black} In Figure~\ref{fig:basis} (bottom, right), we also show the reconstructed data compared to the orginal data using Algorithm \ref{alg:qr} with $M_E=20$ and $M_Q=380$; in this case, the reconstruction error, $\|x-QQ'x\|_2= 5\times10^{-6}$, where columns of $Q\in\mathbb{R}^{1000\times 400}$ are orthogonal.}

Given $M$ sufficiently large, the error in the approximation in \eqref{DF:approxA1} is bounded above by:
\BEA
|\tilde{A}_{kj}-A_{kj}| \leq -\lambda_k b_0 \sqrt{\tau} + \mathcal{O}(\tau),\label{errorbound}
\EEA
where $b_0 \geq \|b(x)\|$ is an upper bound of the diffusion tensor $b(x)$ in \eqref{SDE}, $\lambda_k\leq0$ is the $k$th eigenvalue of the operator $\hat{\cal L}$ with the following order, $0=\lambda_0\geq \lambda_1\geq \ldots$. Note that $\lambda_k\leq0$ is also a minimizer of the Dirichlet norm, $\|\nabla f \|_{\peq}$ for any $f\in H^2(\mathcal{M},\peq)\cap \mathcal{H}_{k-1}^\perp$, where $\mathcal{H}_{k-1}$ denotes the eigen subspace spanned by the first $k-1$ leading eigenfunctions of the operator $\mathcal{\hat L}$ (see the Appendix of \cite{bgh:15} or Chapter 6 of \cite{harlim:18} for the detailed derivations). With the new proposed orthogonal basis, this error bound still holds with a larger positive constant replacing $-\lambda_k\geq 0$ for stochastic problem. In the case of stochastic systems (i.e., $b_0 \neq 0$), numerical examples show that the forecast skill with the new basis functions will not be superior compared to using a large set of eigenfunctions when they are available. However, using the new basis functions to augment existing few eigen functions can improve the forecasting results. For the deterministic problem, the first term in the right hand side of \eqref{errorbound} vanishes since $b_0=0$. Therefore, the error bound is on the order of $\tau$ with any orthogonal basis functions including the newly proposed basis. In the numerical section, we will demonstrate that the forecasting skill based on the new basis functions is superior compared to that using the same number of eigenbasis.

\section{Forecasting dynamical systems}\label{sec4}

In this section, we demonstrate the proposed skill of the diffusion forecasting model in predicting deterministically chaotic and stochastic dynamical systems. In particular, we will consider three examples: a three-dimensional Lorenz model \cite{lorenz:63}, a six-dimensional Lorenz-96 model \cite{lorenz:96}, and a triad stochastic model that is a canonical model for turbulent dynamics \cite{majda:2016}. In the deterministic examples, we will see significant forecast improvements using the new basis functions obtained from the QR decomposition relative to the one that uses purely eigenfunctions. For the stochastic example, similar improvement is found relative to the forecasts using a small set of purely eigenfunctions.

A crucial component for skillful forecasting of chaotic dynamical systems is the accuracy of the initial conditions. For the diffusion forecasting model in \eqref{DF:c}, we need to specify the coefficients: 
\BEA
c_k(0) = \langle p_0, \varphi_k \rangle, \nonumber
\EEA
which is the representation of the initial density $p_0$ in the coordinate basis $\varphi_k(x)$. Practically, however, the initial density, $p_0(x)$, is rarely known. Instead, one usually needs to specify the initial density from some observations $y$ of the current state. Mathematically, one would like to estimate the conditional density $p(x|y)$ and use it as an initial condition for the diffusion forecasting model. 

In the remainder of this section, we will first discuss two different methods to specify the initial densities. One of them is the Bayesian inversion introduced in \cite{bh:16physd} for noisy observations with known error distribution. The second one is a new algorithm based on an application of the Nystr\"om extension \cite{nystrom:30} for clean data set. Subsequently, we discuss the experimental design. We close this section with the numerical results on the three test examples.

\subsection{Specifying initial conditions}
First, let's review the Bayesian inversion technique introduced in \cite{bh:16physd}. Given noisy observations,
\BEA
y_n = x_n + \eta_n, \quad n=1,\ldots,N_V,\nonumber
\EEA
with known error distribution, one can specify the initial density by iterating the following predictor-corrector steps: 
\BEA
p(x_n) &=& e^{\tau\mathcal{L^*}} p(x_{n-1}|y_{n-1}) \label{predictor}\\
p(x_n|y_n) &\propto& p(x_n) p(y_n|x_n).\label{bayes}
\EEA
In \eqref{predictor}, the prior density, $p(x_n)$, at time $n$ is obtained by solving the diffusion forecasting model in \eqref{fouriersum}-\eqref{DF:c} starting from the posterior density at the previous time step $p(x_{n-1}|y_{n-1})$, where the time lag is defined as $\tau = t_n-t_{n-1}$. The Bayes' theorem in \eqref{bayes} is used to specify the posterior density $p(x_n|y_n)$. Since the observation error distribution is known, the likelihood function in \eqref{bayes} is specified as follows $p(y_n|x_n) = p(y_n-x_n)=p(\eta_n)$. In our implementation, we will start this iterative process with a uniform prior $p_1(x)$ and discard the posterior densities at first 10 steps to allow for a spin-up time.  In the Lorenz-63 and the stochastic triad examples below, we will apply this Bayesian inversion with a Gaussian distributed observation error. When the distribution is not known, one can still implement the Bayes' correction in \eqref{bayes} with a non-parametric likelihood function estimator as proposed in \cite{bh:17mwr}.

When the observations are noise-less, then the initial density is given by $p(x_n|y_n) = \delta(x_n-y_n)$, such that the coefficients,
\BEA
c_k(t_n) = \langle p(\cdot |y_n), \varphi_k \rangle = \int_\mathcal{M} \delta(x-y_n) \varphi_k(x) dV(x) = \varphi_k(y_n). \label{initcon}
\EEA
Practically, this requires one to evaluate $\varphi_k$ on a new data point $y_n$ that does not belong to the training data set, $\{x_i\}_{i=1,\ldots, N}$. One way to do this evaluation is with the Nystr\"om extension \cite{nystrom:30}, which is based on the basic theory of reproducing kernel Hilbert space (RKHS) \cite{aronszajn:50}. In our particular application, let $L^2(\mathcal{M},\peq)$ be a RKHS with a symmetric positive kernel $\hat{\cal T}:\mathcal{M}\times\mathcal{M}\to\mathbb{R}$ that is defined such that the components of \eqref{T_hat} are given by $\hat{T}_{ij} = \hat{\cal T}(x_i,x_j)$, for training data $\{x_i\}_{i=1,\ldots,N}$. To be consistent with the conjugation in \eqref{conjugation}, we also define a kernel function ${\cal T}:\mathcal{M}\times\mathcal{M}\to\mathbb{R}$ such that the component of the non-symmetric matrix in \eqref{eqn:T} is given as $T_{ij} = {\cal T}(x_i,x_j) = \hat{q}^{-1/2}(x_i)\hat{\cal T}(x_i,x_j)\hat{q}^{1/2}(x_j)$, where $\hat{q}(x_i) = \frac{1}{N}\sum_{j=1}^N \hat{K}_\epsilon(x_i,x_j)$. Then for any function $f\in L^2(\mathcal{M},\peq)$, the Moore-Aronszajn theorem states that one can evaluate $f$ at $a\in\mathcal{M}$ with the following inner product, $f(a) = \langle f, \hat{\cal T}(a,\cdot)\rangle_{peq}$. 

In our application, this amounts to evaluating,
\BEA
\label{eqn:tt}
\varphi_k(y_n) &=& \hat{q}^{-1/2}(y_n) u_k(y_n)  \nonumber \\ &=& \hat{q}^{-1/2}(y_n) \langle u_k, \hat{\cal T}(y_n,\cdot)\rangle_{peq}\nonumber \\  &\approx&  \hat{q}^{-1/2}(y_n) \frac{1}{N} \sum_{i=1}^N\hat{\cal T}(y_n,x_i) u_k(x_i) \nonumber \\
 &=& \hat{q}^{-1/2}(y_n) \frac{1}{N} \sum_{i=1}^N\hat{\cal T}(y_n,x_i)\hat{q}^{1/2}(x_i) \varphi_k(x_i)\nonumber \\
 &=& \frac{1}{N} \sum_{i=1}^N \mathcal{T}(y_n,x_i)\varphi_k(x_i).
\EEA
In \eqref{eqn:tt}, we have used the fact that $u_j(x_i)$ is the $ij-$th component of the matrix $U$ consisting of the eigenvectors of $\hat{T}$ in \eqref{T_hat}, i.e. $\hat{T}U=U \Lambda$; and the fact that the eigenvectors of $T$ are defined via the conjugation formula $\varphi_k(x_i) = \hat{q}^{-1/2}(x_i) u_k(x_i) $. In summary, the estimation of the coefficients in \eqref{initcon} requires the construction of an $N$-dimensional vector whose $i$th component is $\mathcal{T}(y_n,x_i)$ as in constructing the components of \eqref{eqn:T}. We will implement this extension formulation on the Lorenz-96 example as well as on the MJO modes application below. 

\subsection{Experimental design}

In our numerical experiment below, we consider training the diffusion forecasting model using a data set of length $N$, $\{x_i\}_{i=1,\ldots, N}$, where $x_i:=x(t_i)$ are the solutions of the system of differential equations with a uniform time difference $\tau=t_i=t_{i-1}$. We will verify the forecasting on a separate data set $\{x_n\}_{n=1,\ldots,N_V}$, where $x_n\notin \{x_i\}_{i=1,\ldots, N}$ and the time difference is exactly the same as that in the training data set, that is, $\tau=t_{n}-t_{n-1}$. {\color{black} Numerically, we generate these data set by solving the system of differential equations at times $t_i=1,\ldots,N+N_V$, and define the training data set as the solutions at $ i=1,\ldots,N$ and the verification data as the solutions at $i=N+1,\ldots,N+N_V$.}

To specify initial conditions, we define observations,
\BEA
y_n = x_n + \eta_n, \nonumber
\EEA 
with two configurations.  For the Lorenz-63 and the triad examples below, we let the observation noises be Gaussian with variances 1 and 0.25, respectively. For the Lorenz-96 example, we let the observations be noise-less, i.e., $\eta_n=0$.  

As noted in Section~\ref{sec3}, the diffusion forecasting model in \eqref{DF:matrixA} will be generated using $M=M_E+M_Q$ modes, where $M_E$ denotes the number of eigenvectors and $M_Q$ denotes the number of columns of $T$ that are used as specified in Algorithm~\ref{alg:qr}. $M_Q=0$ means that  only eigenvectors are used as basis functions. On the other hand, the experiments with modes coming from only the columns of $T$ will be denoted with $M_E=0$. The experiments with the mixed modes are denoted by nonzero $M_E$ and $M_Q$.

To measure the forecasting skill, we compute the root-mean-square-error (RMSE) of the mean estimate as a function of forecast lead time $j\tau$,
\BEA
E(j\tau) = \Big(\frac{1}{N_V-10}\sum_{n=11}^{N_V} (\mathbb{E}[x_{n+j}] - x(t_{n+j}) )^2\Big)^{1/2},\label{RMS}
\EEA
where $\mathbb{E}[x_{n+j}] = \int_\mathcal{M} x e^{j\tau\mathcal{L}^*}p(x|y_n) dV(x)$ is the mean forecast statistics computed using the diffusion forecasting model. In \eqref{RMS}, the error is averaged over the verification period, ignoring the first 10 steps to allow for a spin up time in specifying the initial conditions.

As a reference, we also show the ensemble forecasting (or Monte Carlo) skill when the underlying model is known. In all of our simulations below, we will use an ensemble of size $1000$. To specify the initial conditions, we use the Ensemble Transform Kalman Filter \cite{bishop:01,hunt:07} when the observation noises are Gaussian. When the observations are noise-less, we just specify the initial ensemble by corrupting $x_n$ with Gaussian noises with variance 0.04.  With such a small perturbation, we hope to depict the sensitivity of the deterministic dynamical system to initial conditions. In this case, the forecasting statistics are computed by averaging over the ensemble of solutions, $\{x^{(k)}_{n+j}\}_{k=1,\ldots,1000}$.

To quantify the forecast uncertainty, we compare the time evolution of the uncentered second order moment from the diffusion forecasting model and that from the ensemble forecasting using the following RMSE metric,
\BEA
E_2(j\tau) = \Big(\frac{1}{N_V-10}\sum_{n=11}^{N_V} (\mathbb{E}[x^2_{n+j}] - \overline{x^2_{n+j}} )^2\Big)^{1/2},\label{RMS2}
\EEA
where  $\mathbb{E}[x_{n+j}] = \int_\mathcal{M} x^2 e^{j\tau\mathcal{L}^*}p(x|y_n) dV(x)$ is the uncentered second order moment computed using the diffusion forecasting model and $\overline{x^2_{n+j}} = \frac{1}{1000}\sum_{k=1}^{1000} (x^{(k)}_{n+j})^2$ is the empirically estimated ensemble forecast uncentered second order moment.
 
\subsection{Lorenz-63 example}
In this section, we test the diffusion forecasting method with new basis functions on the famous Lorenz-63 model \cite{lorenz:63}. The governing equation of the Lorenz-63 model is given as follows,
\BEA
\frac{dx}{dt} &=& \sigma(y-x) \nonumber \\ 
\frac{dy}{dt} &=& \rho x-y-xz \label{lorenz63} \\ 
\frac{dz}{dt} &=& xy-bz \nonumber
\EEA
with standard parameters $\sigma=10$, $b=8/3$, and $\rho=28$. The dynamics of Lorenz-63 has a chaotic attractor with one positive Lyapunov exponent 0.906 that corresponds to a doubling time of 0.78 time units. Both the training and verification data are generated by solving the ODE in \eqref{lorenz63} with RK4 method with $\Delta t=0.01$. In all of the experiments below, we set $\tau=0.1$. 

\begin{figure}
  \centering
  \includegraphics[width=.8\linewidth]{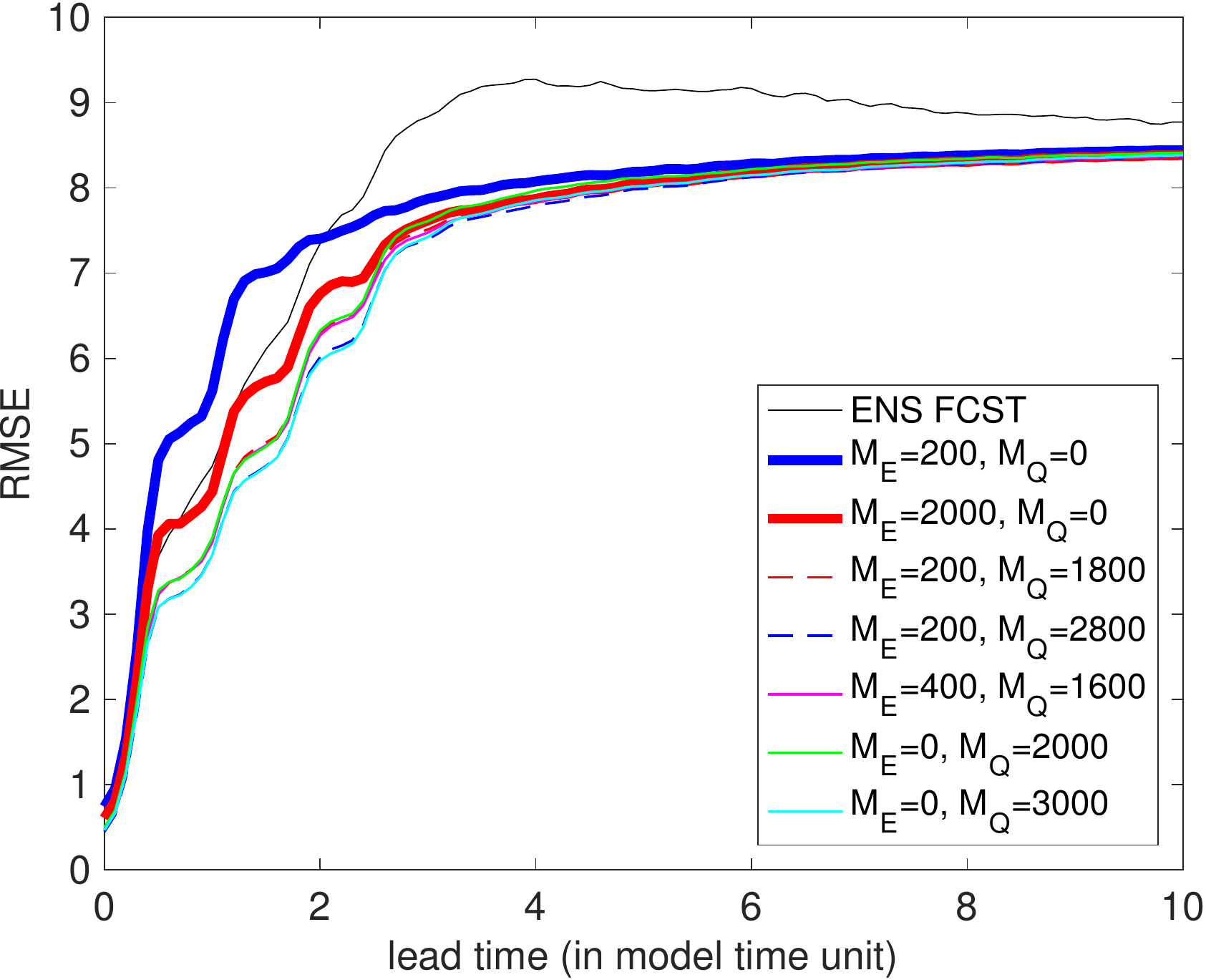}\\
  \includegraphics[width=.8\linewidth]{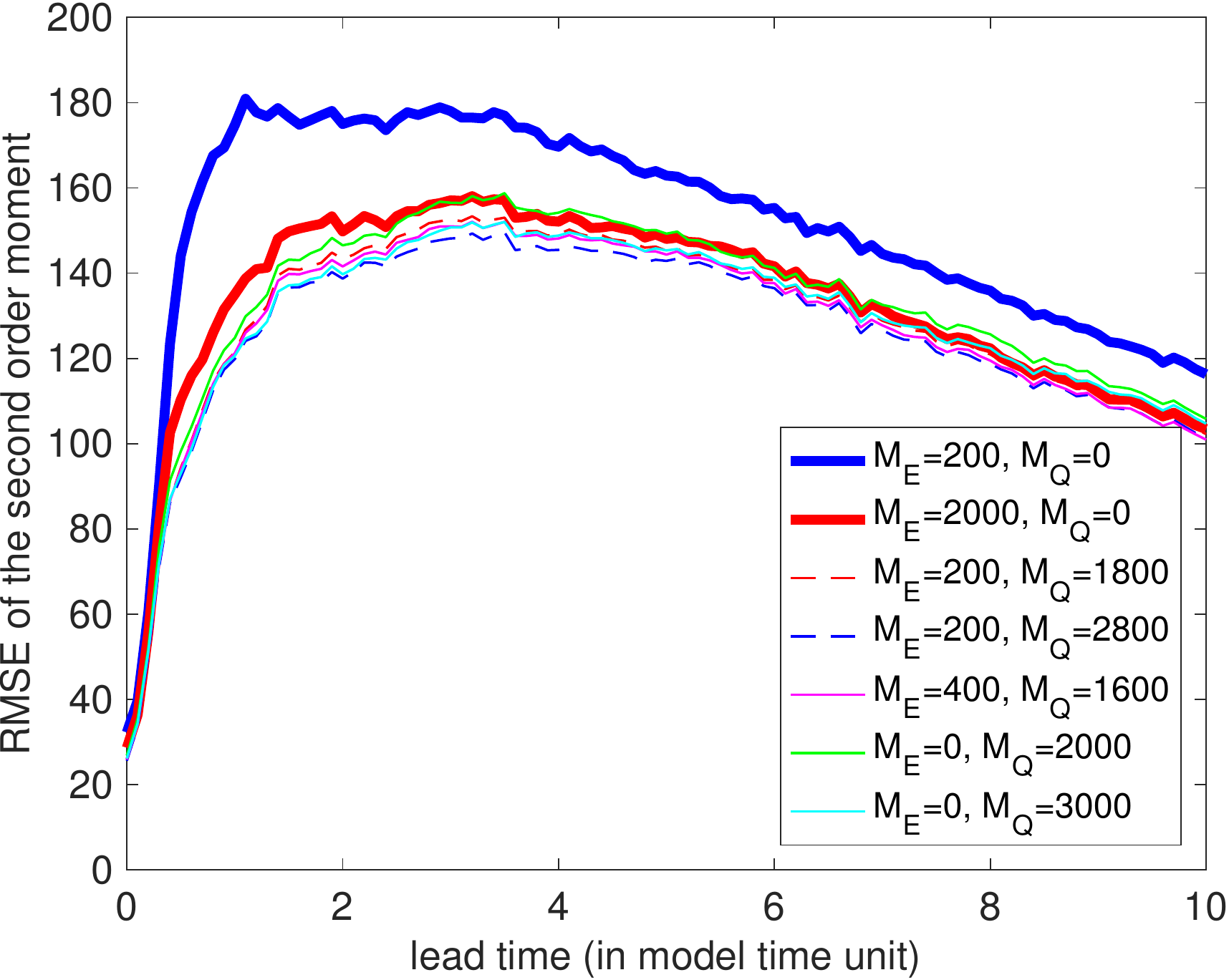}
  \caption{RMSE of the mean (left) and second-order uncentered moment (right) of the diffusion forecasting with various choices of basis functions on the Lorenz-63 example.}
\label{L63fig1} 
\end{figure}

\begin{figure}
  \centering
  \includegraphics[width=.9\linewidth]{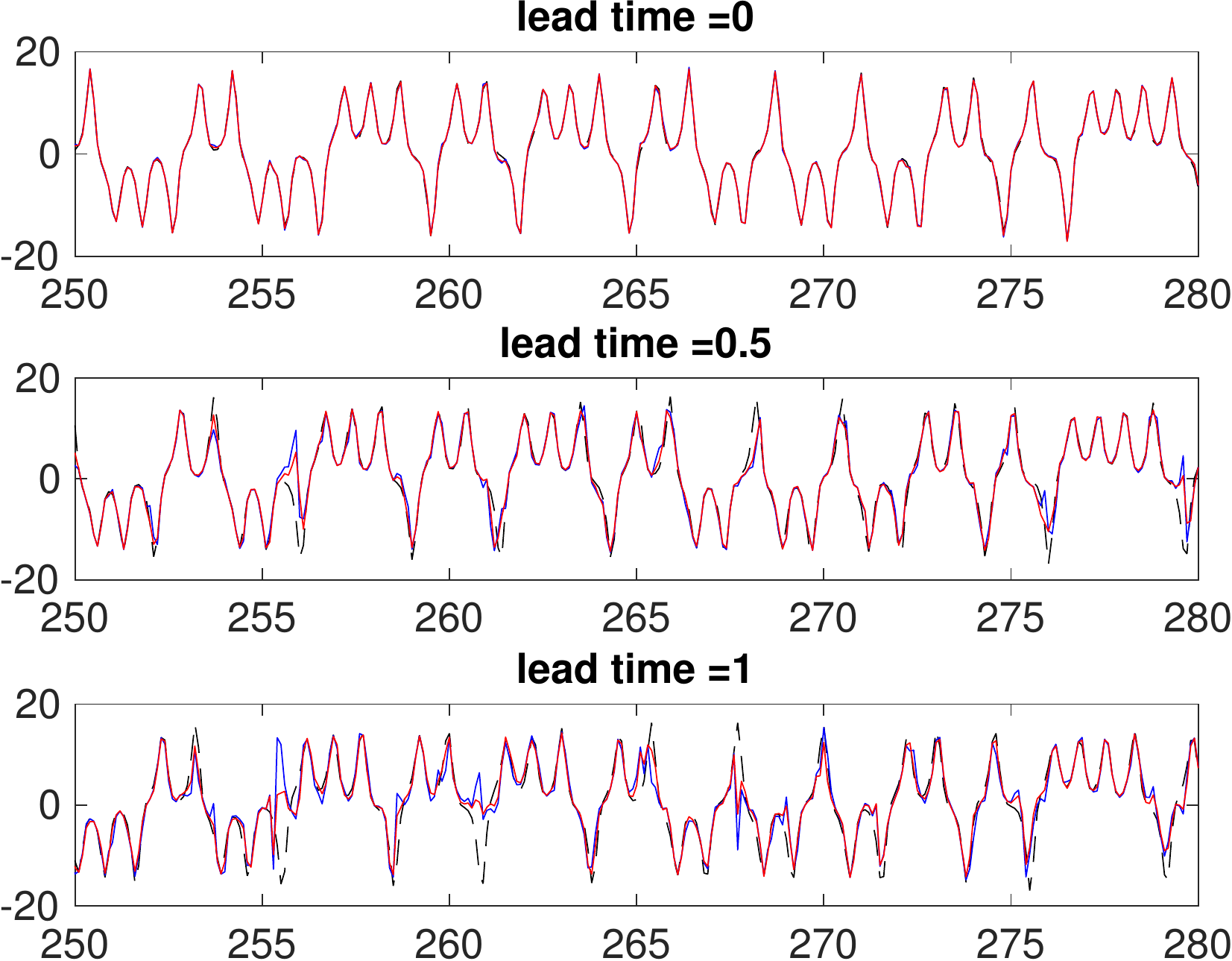}\\
  \includegraphics[width=.75\linewidth]{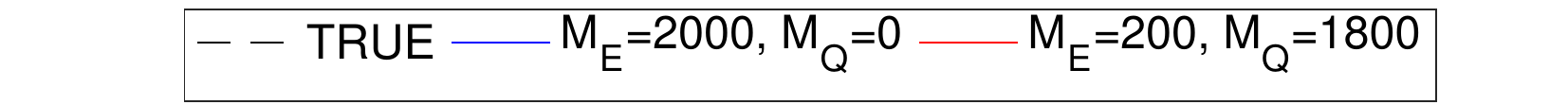}
  \caption{Comparison of prediction on the verification time interval $[250,280]$ with different choices of basis functions. Diffusion forecasting mean estimates of the first component compared to the truth at initial time and forecast lead times 0.5 and 1 model time unit.}
\label{L63fig2}   
\end{figure}

\begin{figure}
  \centering
  \includegraphics[width=.9\linewidth]{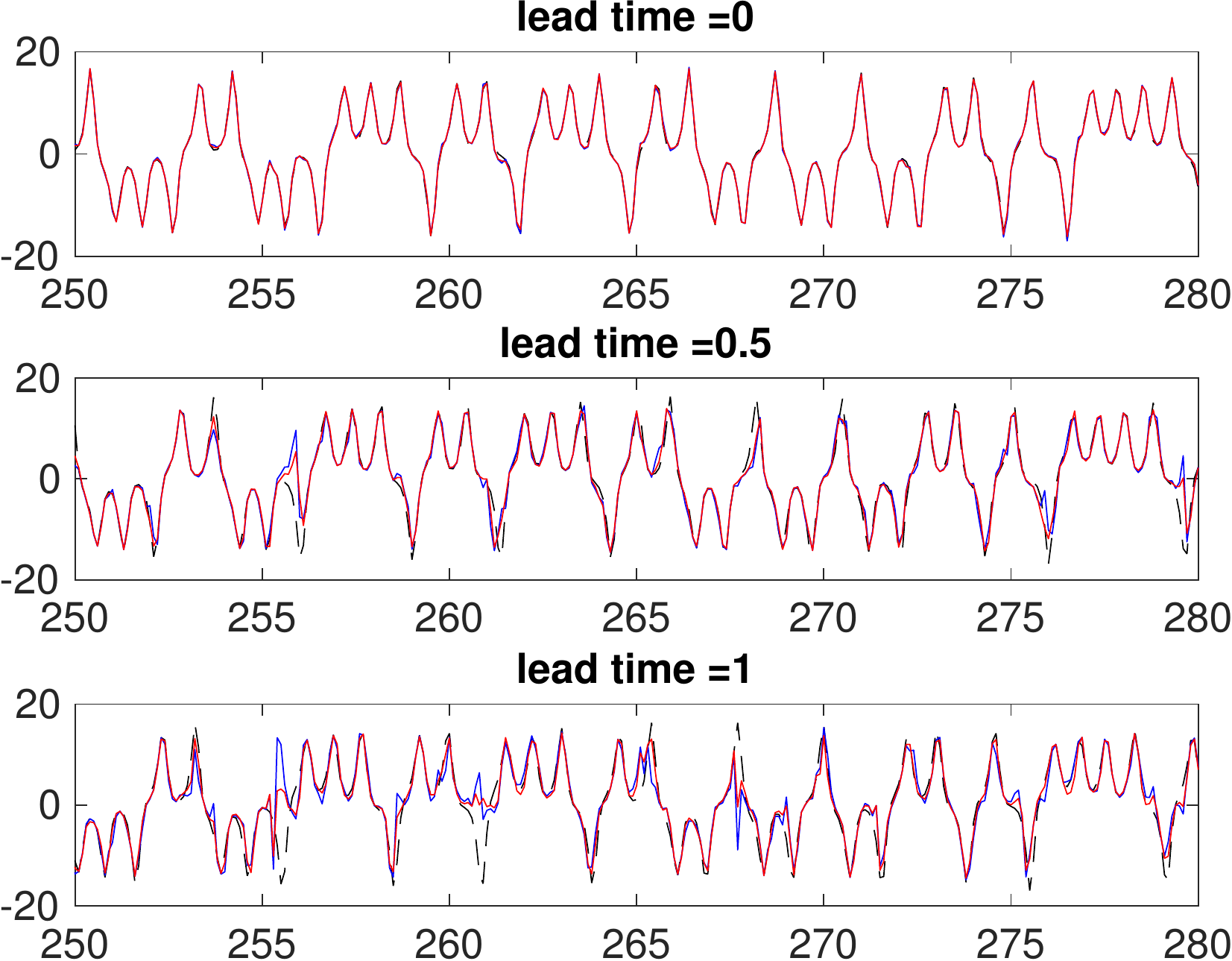}\\
  \includegraphics[width=.75\linewidth]{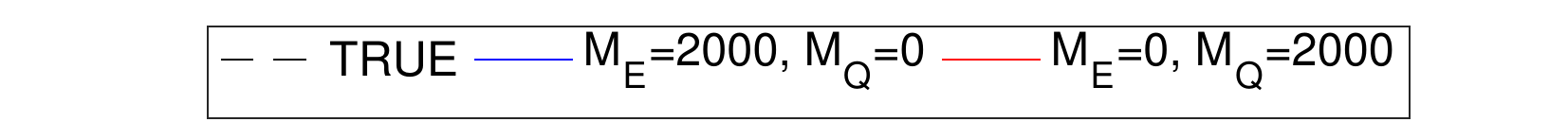}
  \caption{Same as Figure~\ref{L63fig2}, except we show results with the purely QR basis, $M_E=0, M_Q=2,000$.}
\label{L63fig2_2}   
\end{figure}

In this example, the training data length is $N=25,000$ and the verification data length is $N_V=10,000$. In Figure~\ref{L63fig1}, we show the RMSEs of the first and second order moment estimates that are computed based on \eqref{RMS} and \eqref{RMS2}, respectively. There are mainly two sets of results in Figure \ref{L63fig1}: one for $M=M_E+M_Q=2000$; another one for $M=3000$. The set for $M=2000$ shows that: the proposed basis functions by QR decomposition in Algorithm \ref{alg:qr} improve the forecasting results relative to those that purely using eigenvectors based on these two metric evaluations; increasing $M_E$ when $M$ is fixed does not improve the forecasting results, i.e., merely using the basis functions from the columns of $T$ in Equation \eqref{eqn:T} are sufficient to make a good prediction when $M_Q$ is large enough. The comparison of the results when $M=2000$ and $M=3000$ shows that increasing $M_Q$ can improve the forecasting accuracy. It is worth emphasizing that the proposed basis leads to better forecasting results than the ensemble forecasting method at all times, while the one of purely eigenfunctions is not.

In Figure~\ref{L63fig2}, we show the estimates at $300$ verification data points on interval $[250,280]$ for the forecasts that use purely eigenvectors $M_E=2000, M_Q=0$ and mixed modes $M_E=200$, $M_Q=1800$ at the initial time and forecast lead-time of 0.5 and 1. {\color{black}In Figure ~\ref{L63fig2_2}, we shows the estimates for the forecasts that use purely eigenvectors $M_E=2000, M_Q=0$ and mixed modes $M_E=0$, $M_Q=2000$. }Notice that these times are before and beyond the doubling time of the forecast initial errors, respectively. At the initial time (lead time-0), notice that the Bayesian inversion in \eqref{predictor}-\eqref{bayes} produces very accurate initial conditions. At lead time-0.5, the mixed-mode forecasts {\color{black}and even the purely QR-mode forecasts} are closer to the truth relative to that of the eigenvectors. The differences between the two are more apparent at lead time-1. {\color{black}In fact, careful inspection shows that the purely QR-mode forecasts are almost the same as the mixed-mode forecasts. Hence, in this example, the purely QR-modes are good enough and it is not necessary to compute the leading eigenfunctions.}

Different tests with varying training data sizes (from $5,000$ to $50,000$) have been performed using the same comparison as discussed just above. The same conclusions above can be drawn from these experiments and results are not shown. 

\subsection{Lorenz-96 example}
In this section, we show results for the Lorenz-96 model \cite{lorenz:96}. The governing equation of the Lorenz-96 model is given as follows,
\BEA
\frac{d x_j}{dt}=(x_{j+1}-x_{j-2}) x_{j-1} -  x_j + F,\quad j=1,\ldots, d,\label{lorenz96}
\EEA
where we set $d=6$ and $F=8$ in the experiment here. Note that this model is designed to satisfy three basic properties: it has a linear dissipation (the $-{x}_j$ term) that decreases the total energy defined as $E=\frac{1}{2}\sum_{j=1}^{d}{x}_j^2$, an external forcing term $F>0$ that can increase or decrease the total energy, and a quadratic discrete advection-like term that conserves the total energy. In our numerical experiment below, both the training and verification data are generated by solving the ODE in \eqref{lorenz96} with RK4 method with time step $\Delta t=0.05$ and we set $\tau=0.05$. 

\begin{figure}
  \centering
  \includegraphics[width=.45\linewidth]{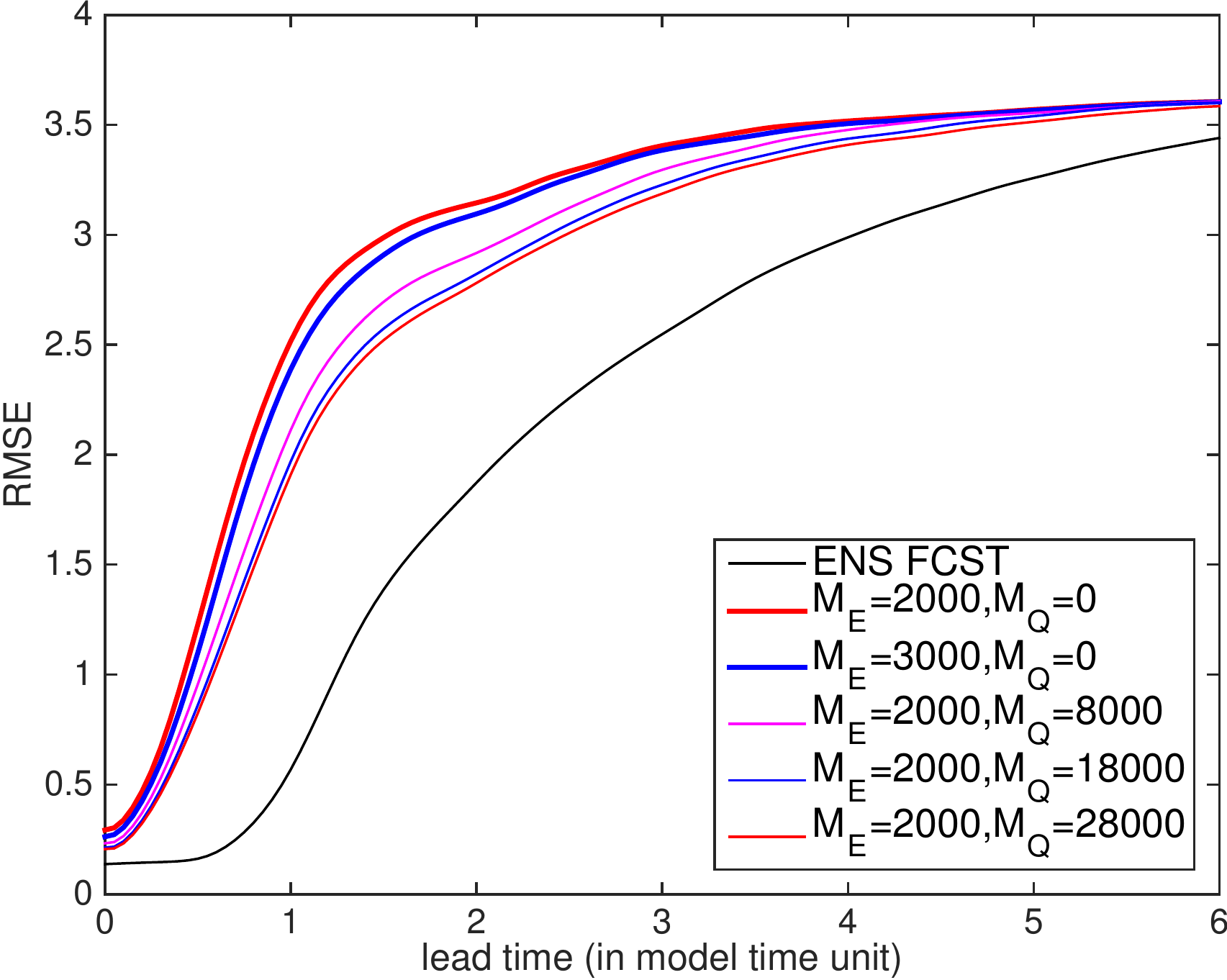}\quad\quad
  \includegraphics[width=.45\linewidth]{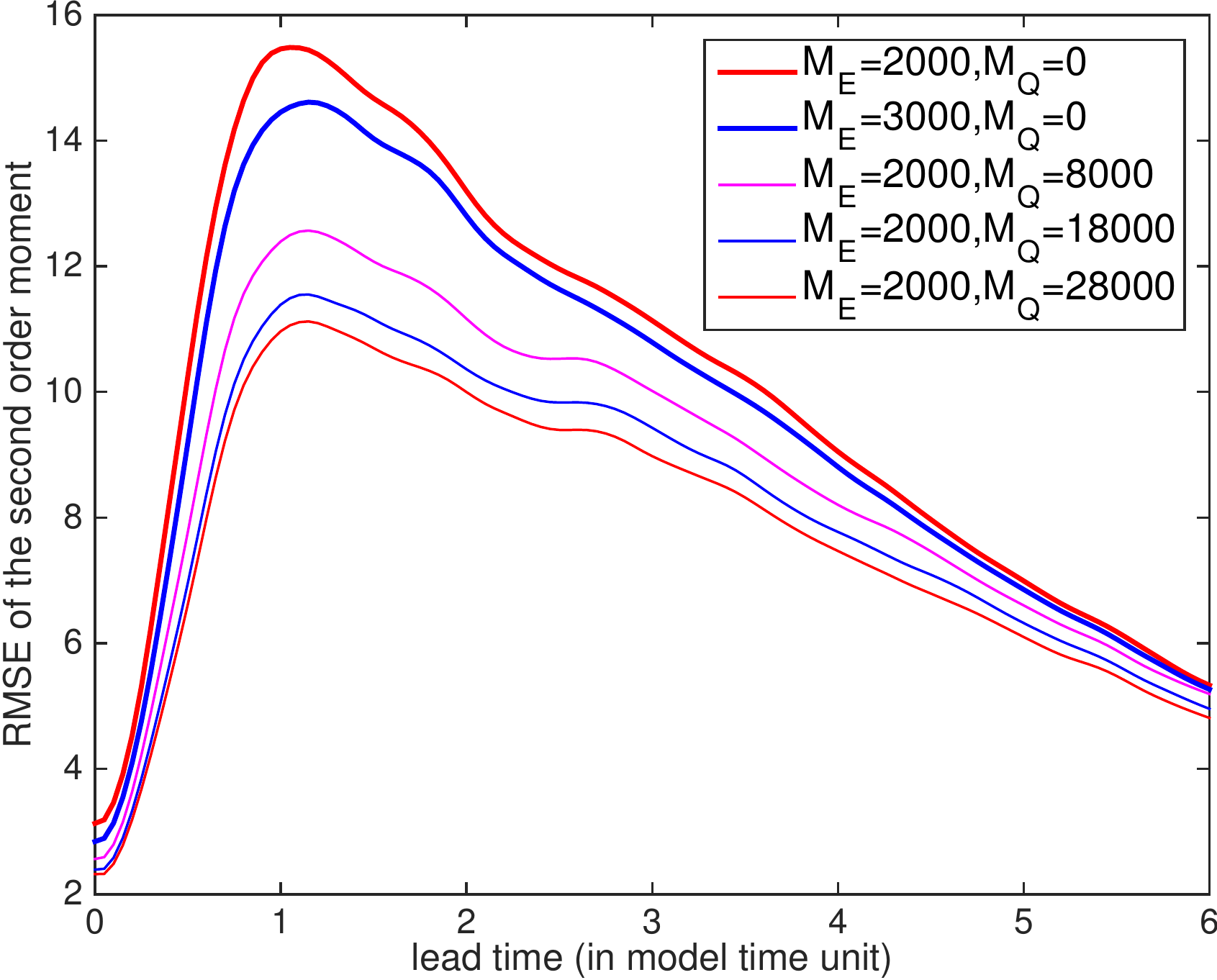}
  \includegraphics[width=.45\linewidth]{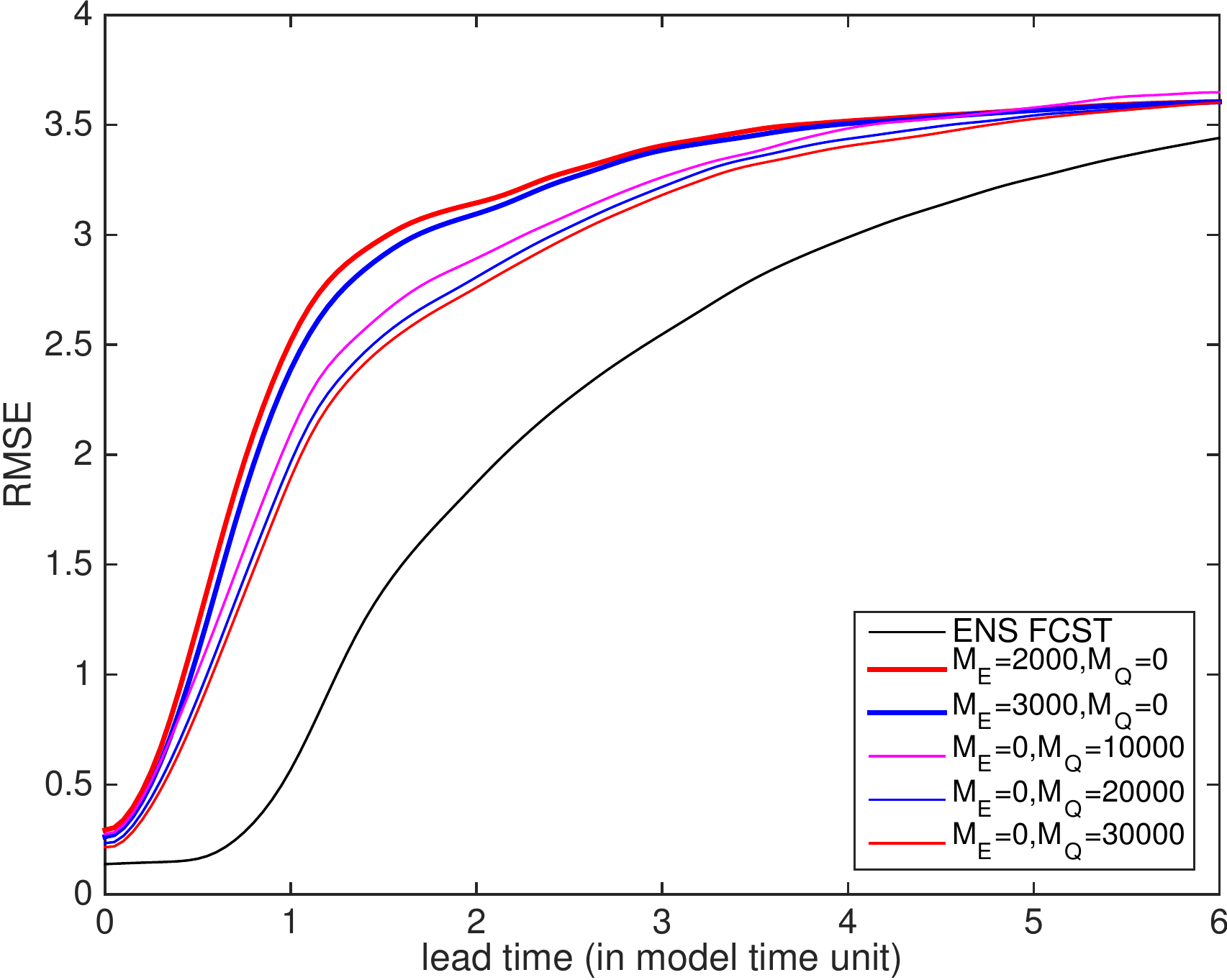}\quad\quad
  \includegraphics[width=.45\linewidth]{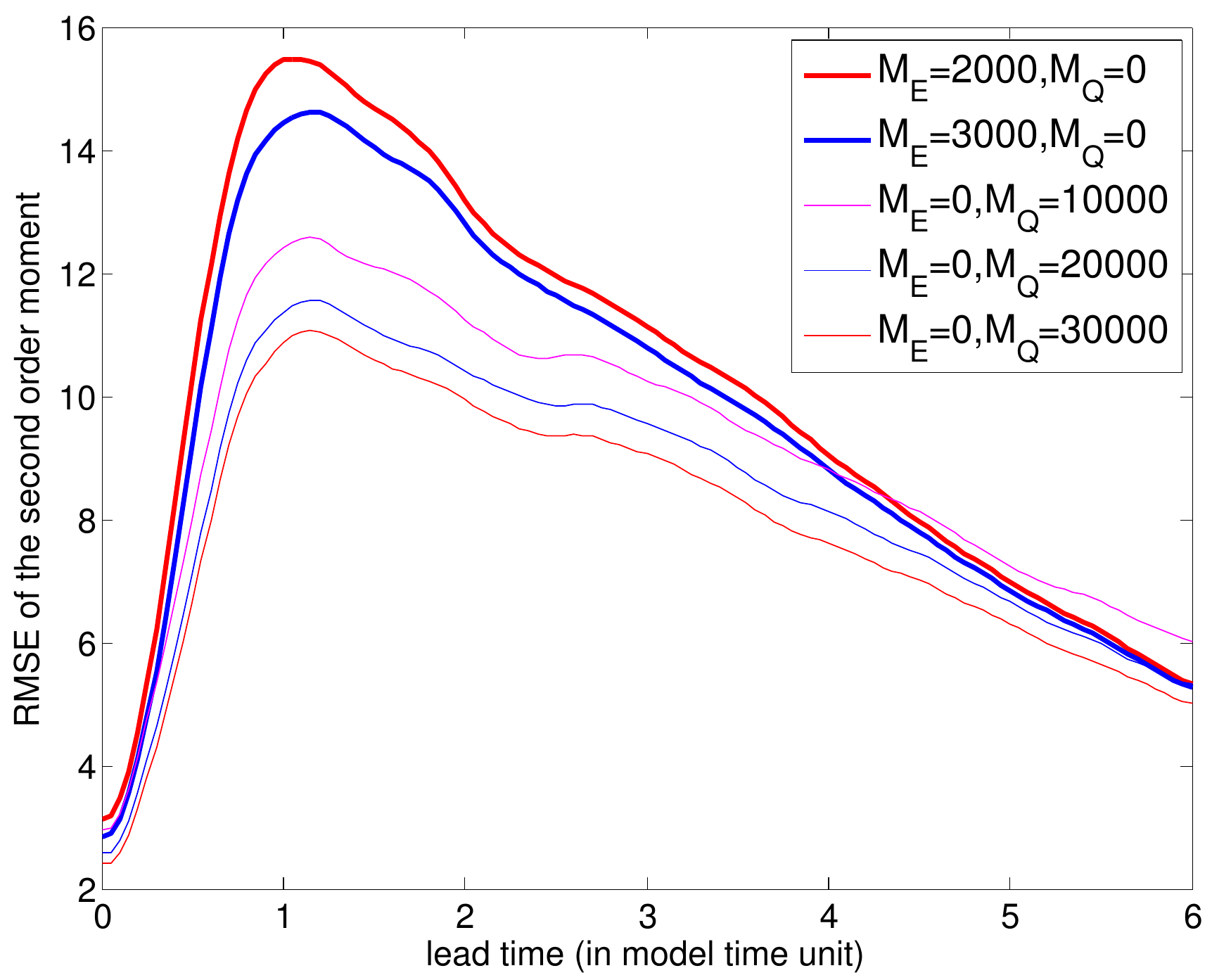}
  \caption{RMSE of the mean (left) and second-order uncentered moment (right) of the diffusion forecasting with various choices of basis functions on the Lorenz-96 example. Top panels show mixed basis, the bottom panels show purely QR basis. Notice the slightly worse with the purely QR basis when the number of total basis functions are 10,000.}
\label{L96fig1} 
\end{figure}

Here, we show results with training data of length $N=100,000$ and verification data of length $N_V=10,000$. For this size of data set, the computing time for the basis functions is recorded in Table~\ref{L96tab1}. Notice that the computational cost for the leading eigenvectors is significantly more expensive than that for QR-decomposition of tall matrices. In Figure~\ref{L96fig1}, we show the RMSEs of the first and second order moment estimates that are computed based on \eqref{RMS} and \eqref{RMS2}, respectively. Notice that the proposed mixed basis functions improve the forecasting results relative to those purely based on eigenvectors according to these two metric evaluations. We should point out that the results of mixed basis with $M_E=2000$ and $M_E=3000$ are indistinguishable. So we only report the results of mixed basis with $M_E=2000$.
{\color{black}The results with the mixed basis are also mostly indistinguishable from those obtained using purely QR basis functions with the same number of total modes (see the bottom part of Figure~\ref{L96fig1}); the purely QR basis models produce larger error at longer time when the total number of basis functions are 10,000 and 20,000.}

\begin{figure}
  \centering
  \includegraphics[width=.9\linewidth,height=11cm]{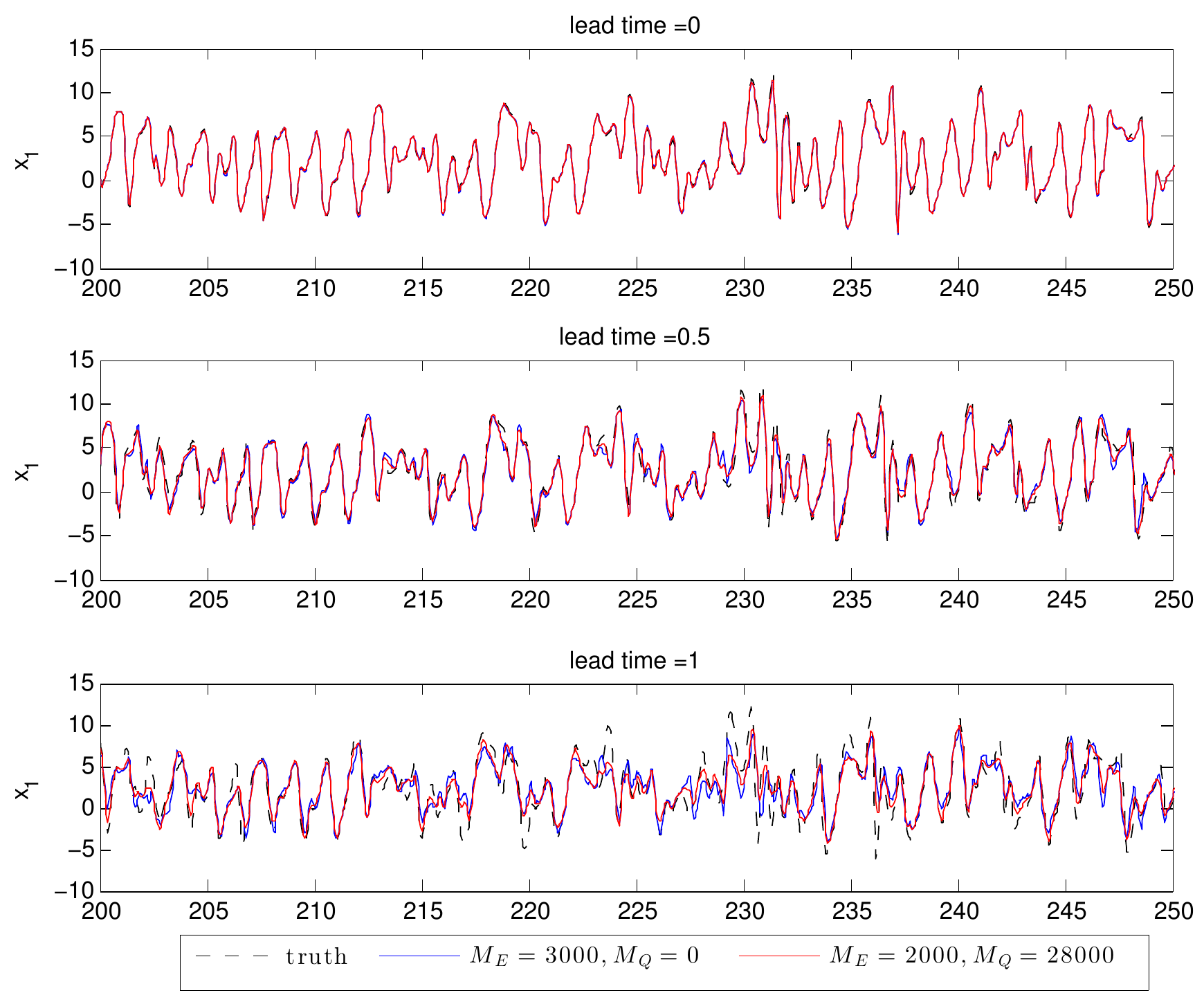}
  \caption{Prediction of the first component of the Lorenz-96 model on the verification time interval $[200,250]$ at initial time and forecast lead times 0.5 and 1 model time unit.}
\label{L96fig2}   
\end{figure}

\begin{figure}
\centering
\includegraphics[width=.9\linewidth,height=11cm]{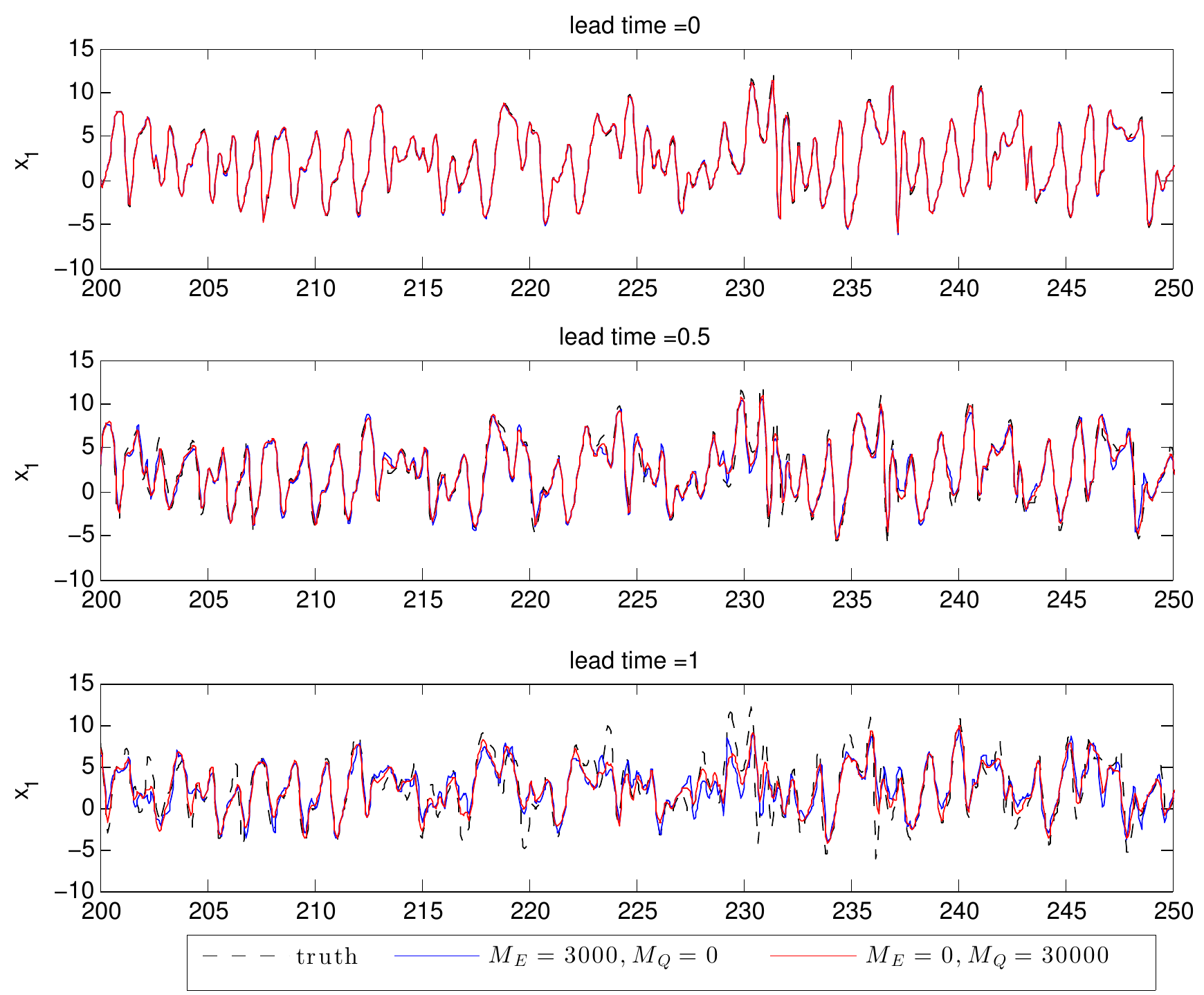}
\caption{Same as Figure~\ref{L96fig2}, except we show results with the purely QR basis, $M_E=0, M_Q=30,000$.}
\label{L96fig2b}   
\end{figure}

Also, the forecasting skill is slightly improved when $M_Q$ in the mixed-mode is increased. In Figure~\ref{L96fig2}, we show the estimates at $1000$ verification data points on interval $[250,300]$ for the forecasts that use purely eigenvectors ($M_E=3000, M_Q=0$) and mixed modes ($M_E=2000$, $M_Q=28000$) at the initial time and forecast lead-time of 0.5 and 1. {\color{black}In Figure~\ref{L96fig2b}, we show the corresponding forecasts using the purely QR modes ($M_E=0$, $M_Q=30000$).} Notice that these times are beyond the doubling time of the forecast initial errors. At the initial time (lead time-0), the Nystrom extension produces very accurate initial conditions. At lead time-0.5, the mixed-mode forecasts (which are almost indistinguishable to those of the purely QR modes) are closer to the truth relative to that of the eigenvectors. The differences between the two are more apparent at lead time-1.

{\color{black}Computationally, the mixed basis $M_E=2000, M_Q=28000$ are obtained in roughly 21 hours (wall clock time), adding up 5.5 and 15.7 hours (see Table~\ref{L96tab1}), which is about the same amount of time it takes to get $M_E=3000$ eigen basis; yet the mixed basis representation produces an improved forecast skill. In fact this forecasting skill is almost identical to that obtained using the purely QR basis (5.5 hours). This suggests that for this example, it is not necessary to use eigenvectors. In general, considering that the QR factorization of an $N\times M_Q$ matrix with $M_Q<N$ is computationally much cheaper than computing the $M_Q$ leading eigenvectors of an $N\times N$ matrix (see Table~\ref{L96tab1} for an example of wall-clock time comparison), generating basis functions for diffusion forecasting by Algorithm \ref{alg:qr} using a large enough $M_Q$ is better than the standard diffusion forecasting approach by computing eigenvectors.}

Note that similar trends as above are found in the case of shorter training data sets (namely $N=50,000$ and $25,000$). The only difference is that the forecast errors are larger; {\color{black} the relative forecast error of using $N=50,000$ compared to $N=100,000$ under the same configuration ($M_E=2000, M_Q=28000$) is about 20\% worse at initial times and gradually decays to about 5\% as the forecast lead time closes to 1 unit.} Based on these findings, our conjecture is that we may need much longer training data to get results closer to the ensemble forecast in this higher dimensional example.

\begin{table}[htp]
\caption{The computational cost (in wall clock time) of computing basis functions based on a data set of length $N=100,000$. Each case is computed in MATLAB using a single-core processor of Intel Xeon E7-4830 v2 2.2GHz with a 1Tb of RAM.}
\begin{center}
\begin{tabular}{|c|c|c|}\hline
method & modes & wall clock time (in hours) \\ \hline
 & 10,000 modes & 0.7 \\
QR & 20,000 modes & 2.6 \\
 & 30,000 modes & 5.5 \\ \hline
Eig & 2000 modes & 15.7 \\
 & 3000 modes & 21.8 \\ \hline
\end{tabular}
\end{center}
\label{L96tab1}
\end{table}%

\subsection{Stochastic triad example}
In this section, we consider an application of the diffusion forecasting method to the following triad model which was introduced in Chapter 2.3 of \cite{majda:2016}  as a simplified model for geophysical turbulence,
\BEA
\frac{dx}{dt} = B(x,x) + Lx -d\Lambda x + \sigma\Lambda^{1/2} \dot{W},\label{triad}
\EEA
where $x\in \mathbb{R}^3$. The term $B(x,x)$ is bilinear, satisfying the Liouville property, $div_x(B(x,x)) = 0$ and the energy conservation, $x^\top B(x,x) = 0$. The linear operator $L$ is skew symmetric, $x^\top L x = 0$, representing the $\beta$-effect of Earth's curvature. In addition, $\Lambda>0$ is a fixed positive-definite matrix and the operator $-d\Lambda x$ models the dissipative mechanism, with a scalar constant $d>0$. The white noise stochastic forcing term $\sigma\Lambda^{1/2}\dot{W}$ with scalar $\sigma^2>0$ represents the interaction of the unresolved scales.

In the following numerical test, we consider $B(u,u) = (B_1 u_2u_3,B_2 u_1u_3,B_3 u_1u_2)^\top$ 
with coefficients chosen to satisfy $B_1+B_2+B_3=0$. In our numerical experiment below, we set $B_1=.5, B_2=1, B_3=-1.5$, 
\BEA
L = \begin{pmatrix}0 & 1 & 0 \\ -1 & 0 & -1 \\ 0 & 1 & 0 \end{pmatrix}, \quad\Lambda = \begin{pmatrix} 1 & 1/2 & 1/4 \\  1/2 & 1 & 1/2 \\ 1/4 & 1/2 & 1\end{pmatrix},\nonumber
\EEA
and $d=1/2$, $\sigma=1/5$. Here, the truth is generated by Euler-Maruyama discretization with $\Delta t=0.01$.

\begin{figure}
  \centering
  \includegraphics[width=.6\linewidth]{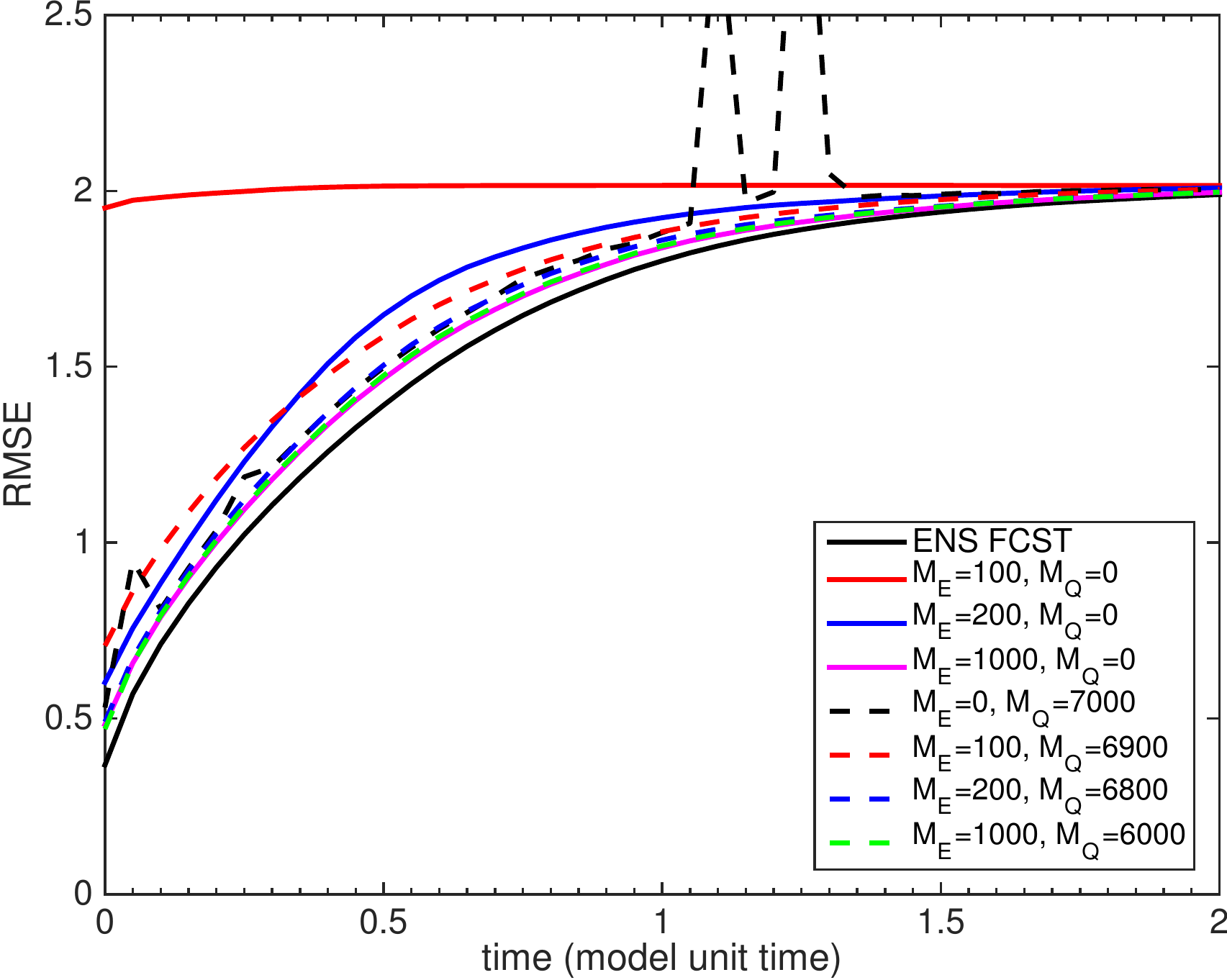}\\
  \includegraphics[width=.6\linewidth]{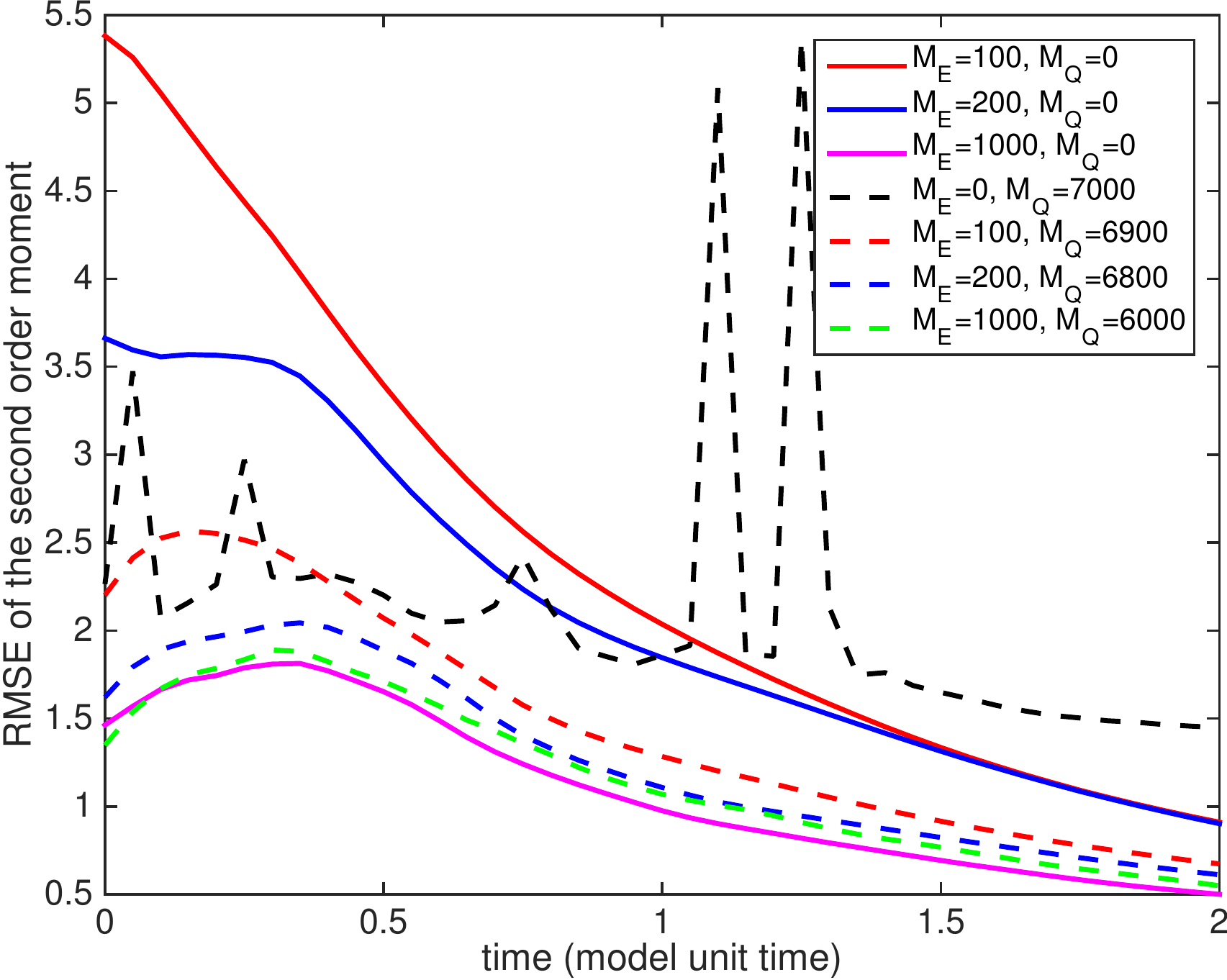}
  \caption{RMSE of the mean (top) and second-order uncentered moment (bottom) of the diffusion forecasting with various choices of basis functions on the triad example.}
\label{triadfig} 
\end{figure}

\begin{figure}
  \centering
  \includegraphics[width=.9\linewidth]{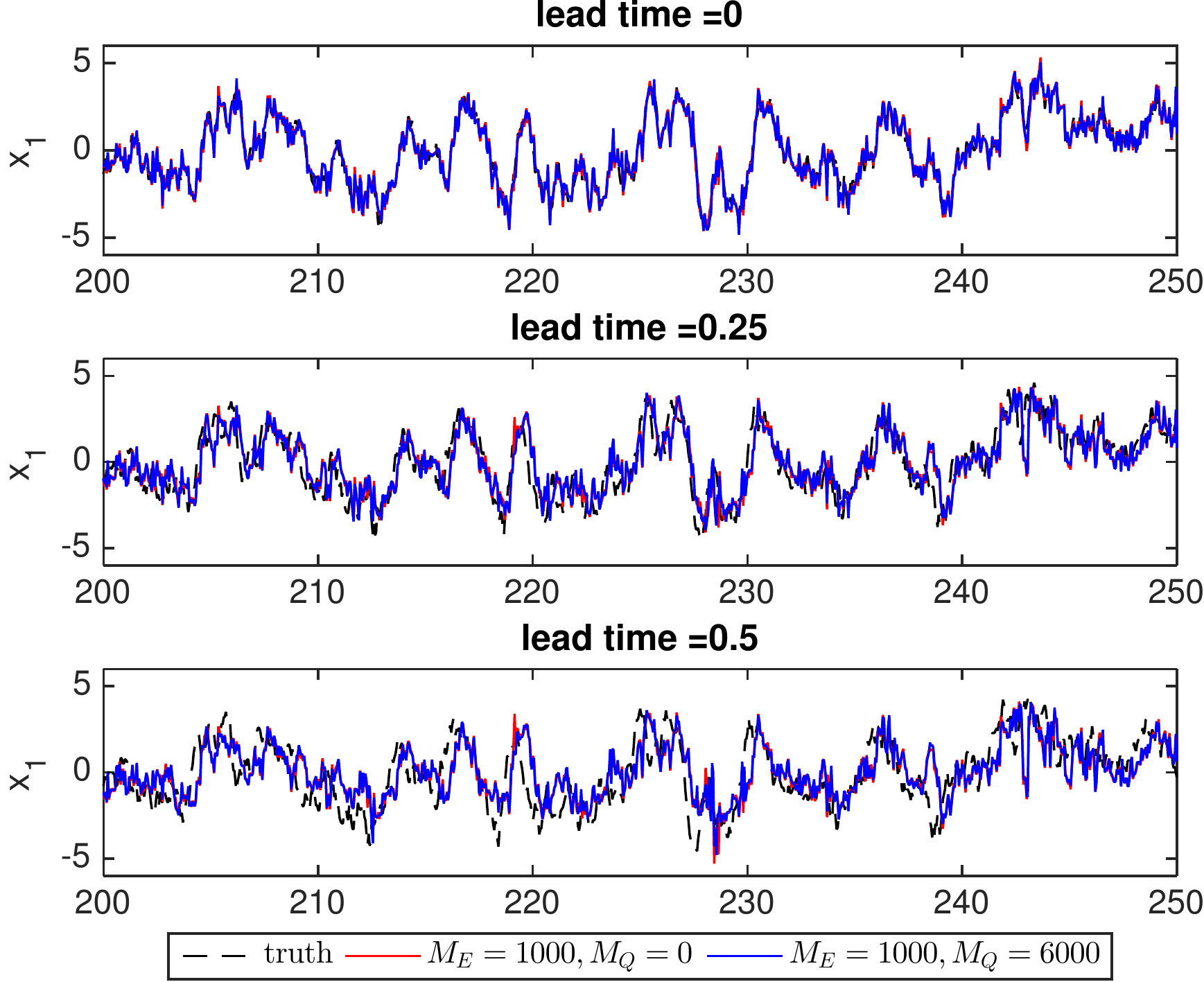}
 \caption{Prediction of the first component of the triad model on the verification time interval $[200,250]$  at initial time and forecast lead times 0.25 and 0.5 model time unit.}
\label{triadtraj} 
\end{figure}

Here, we show results with training data of length $N=25,000$ and verification data of length $N_V=10,000$ with temporal step $\tau=0.05$. In Figure~\ref{triadfig}, we show the RMSEs of the first and second order moment estimates that are computed based on \eqref{RMS} and \eqref{RMS2}, respectively. {\color{black}In Figure~\ref{triadtraj}, we show the corresponding prediction of the first component $x_1$ at the verification interval $[200,250]$ at initial times and forecast lead times 0.25 and 0.5.} First, notice that using purely eigenvectors ($M_E=1000, M_Q=0$) produces forecasting skill close to that of the ensemble forecasting in this example. When the number of eigenvectors is small ($M_E=100$ or $M_E=200$), the forecast skills deteriorate. In the case of deterioration, we can improve the forecast by using mixed modes (see the results for $M=M_E+M_Q=7000$). On the other hand, the forecast skill using purely $M_Q=7000$ QR modes is not stable. However, forecasting with QR modes can be stabilized by using more $M_Q$; but even so, the skill is not superior compared to that with purely $M_E=1000$ eigenvectors. 

The numerical results here confirm the theoretical error estimate in Chapter~6 of \cite{harlim:18}. If the number of modes $M$ is large enough, the $M$-term approximation error in \eqref{fouriersum} is very small; in this case, the dominant term in the forecasting error is the approximation error of \eqref{shiftapprox}, the bound of which is $-\lambda_k b_0 \sqrt{\tau} + \mathcal{O}(\tau)$ given by \eqref{errorbound} when eigenbasis is used in \eqref{DF:approxA1} and \eqref{DF:approxA} (see Theorem $6.2$ in Chapter~6 of \cite{harlim:18} for details). Following the same proof as in Theorem $6.2$ in Chapter~6 of \cite{harlim:18}, if another orthonormal basis is used in \eqref{DF:approxA1} and \eqref{DF:approxA}, the forecasting error is dominated by $cb_0 \sqrt{\tau} + \mathcal{O}(\tau)$ since the Dirichlet norm of the new basis functions is larger than the minimizer, $\|\nabla \varphi_k\|_{\peq}=c\geq -\lambda_k$. Hence, in the case of a large $M$, diffusion forecasting with eigenbasis is the optimal choice for stochastic problems. However, if $M$ is small, i.e., only few eigen functions are available, the $M$-term approximation error in \eqref{fouriersum} is large and becomes the dominant term in the forecasting error. In this case, one can reduce the $M$-term approximation error with the mixed modes by Algorithm \ref{alg:qr} to reduce the final forecasting error. 

\section{Forecasting the MJO spatial NLSA modes}\label{sec5}

In this section, we consider predicting two spatial modes associated with the boreal winter Madden-Julian Oscillation (MJO). These two spatial modes were obtained by applying the nonlinear Laplacian Spectral Analysis (NLSA) algorithm \cite{gtm:12} on the Cloud Archive User Service (CLAUS) version 4.7 multi-satellite infrared brightness temperature and are associated with the MJO modes based on the broad peak of their frequency spectra around 60 days and the spatial pattern that resembles the MJO convective systems around the Indian Ocean western Pacific Sector (see \cite{gtm:12} for the detailed analysis). 

In \cite{cmg:14}, a four-dimensional parametric stochastic model has been proposed to predict this time series extensively. Here, we mimic the study in \cite{cmg:14} with a slightly different configuration. In particular, we will train our model using the NLSA modes from 3-Sep-1983 00:00:00 to 15-Mar-2004 18:00:00, with time discretization $\tau=6$ hours, so the total number of training data points is $N=30,000$. The verification data is the NLSA modes (obtained separately from the training data) from 01-Jul-2006 00:00:00 to 28-Dec-2008 18:00:00, also with time discretization of $\tau=6$ hours and the total number of verification data points is $N_V=3,647$ (both the training and verification data were provided to us by D. Giannakis, who is one of the inventors of the NLSA algorithm). 

\begin{figure}
\centering
\includegraphics[width=.8\linewidth]{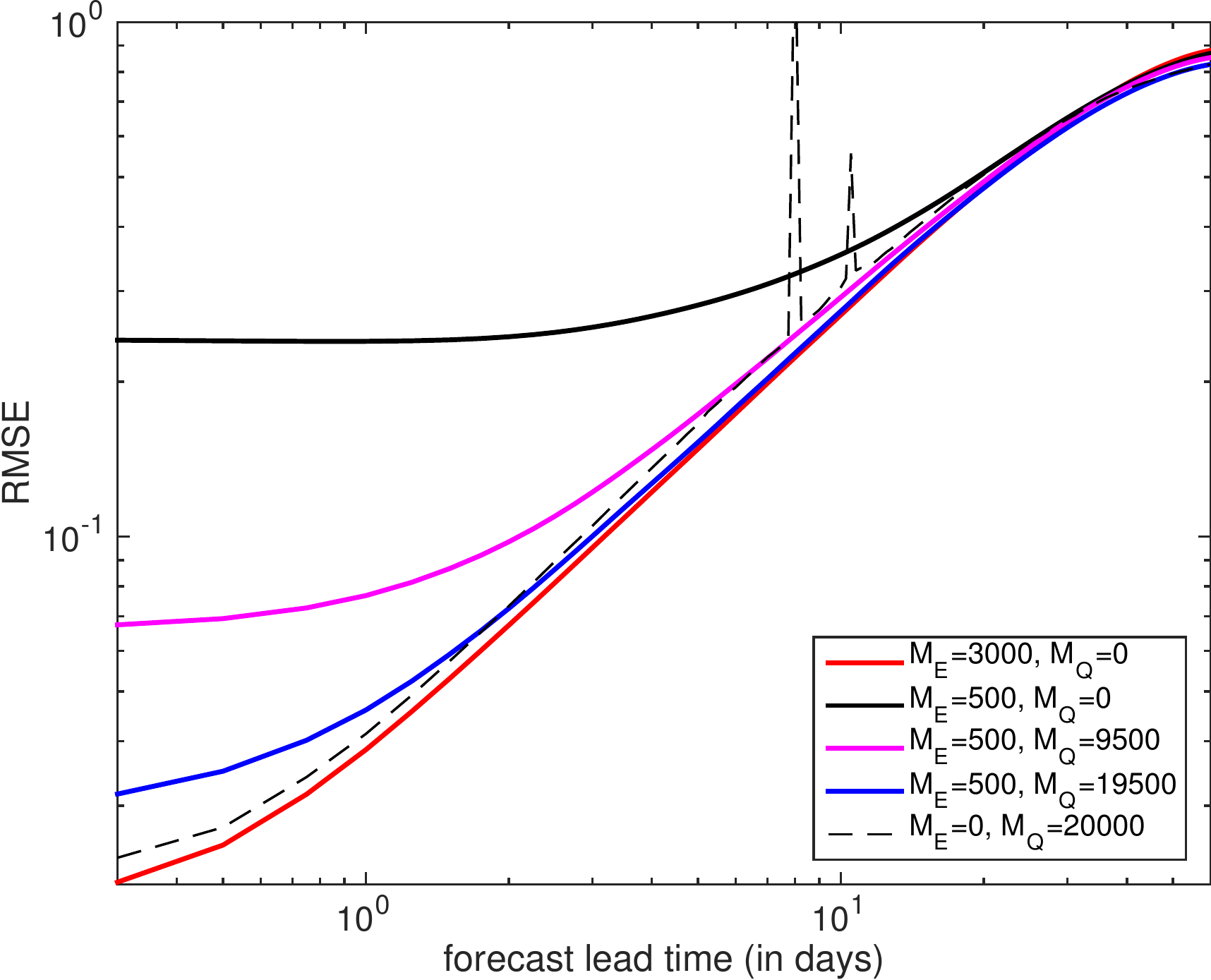}
\caption{RMSE of the mean forecast of the diffusion forecasting with various choices of basis functions on the MJO example.}
\label{MJOfig1} 
\end{figure}

\begin{figure}
\centering
\includegraphics[width=.8\linewidth]{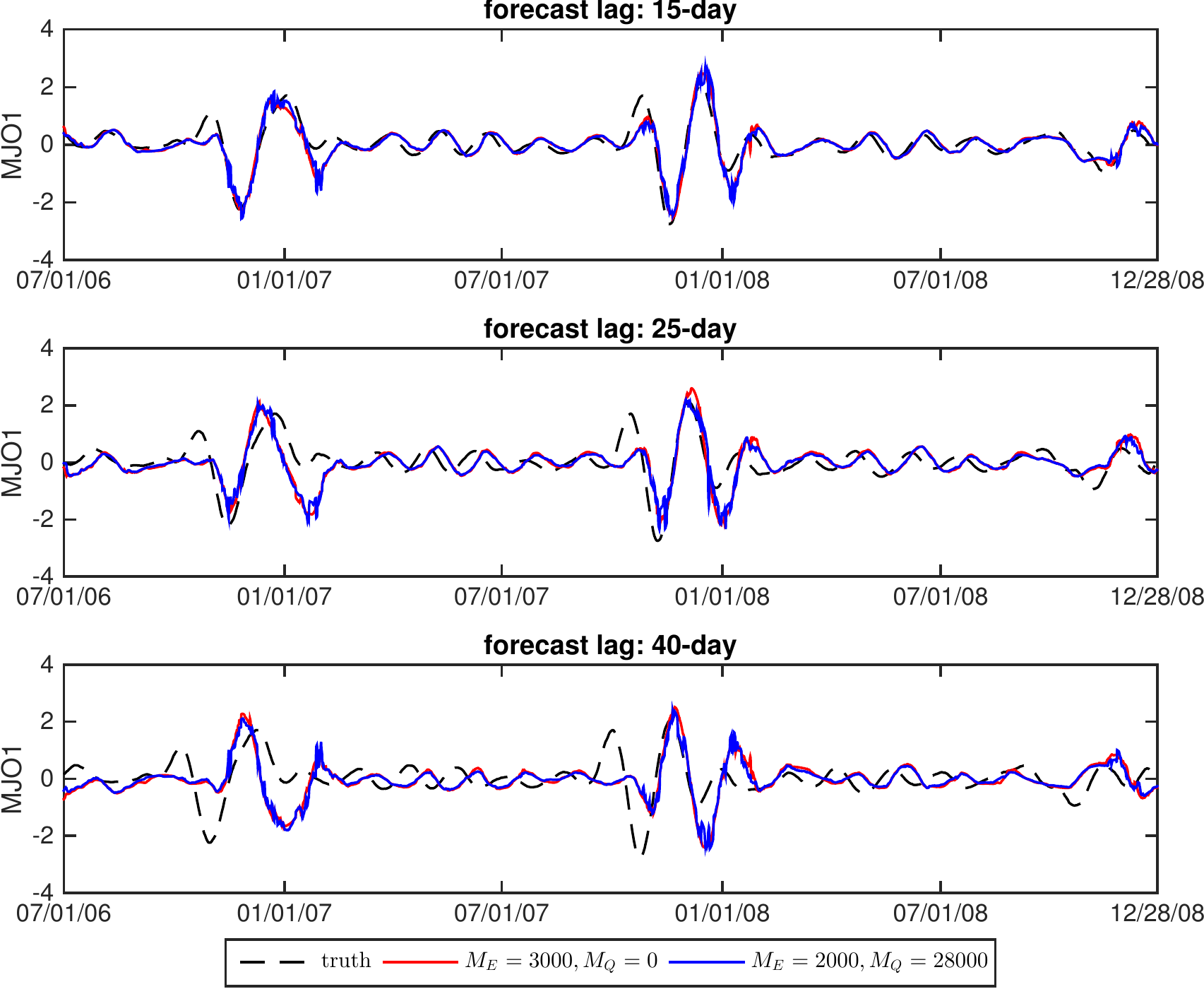}
\caption{Diffusion forecasting mean estimates of the MJO1 compared to the truth at forecast lead times of 15-days, 25-days, and 40-days.}
\label{MJOfig2} 
\end{figure}

In Figure~\ref{MJOfig1}, we show the RMSE of the mean forecast estimates using various basis functions. Notice that the mixed basis with $M_E=500$, $M_Q=9500$ improves the estimates relative to just using purely $M_E=500$ eigenvectors. The forecasting estimates can be further improved by increasing $M_Q$, say $M_E=500, M_Q=19500$. However, when we use only purely $M_Q=20,000$ QR basis, the prediction skill is unstable (dashes). This justifies the use of a few eigen modes to stabilize the performance of the QR basis. In terms of RMSE, notice that the result of $M_E=3000$ is better than that of $M_E=500, M_Q=19500$; we speculate that this behavior is because the underlying dynamic is stochastic. However, these two sets of basis functions produce similar results for trajectories at specific times (5, 25, and 40 days) as shown in Figure~\ref{MJOfig2}). Notice that the forecasting skill for 15 days is quite high and it deteriorates as the lead time increases to 25 and 40 days. 


It is worth pointing out that in our implementation, it is implicitly assumed that these two spatial modes can be modeled with a two-dimensional autonomous dynamic, as opposed to the approach in \cite{cmg:14} that uses a four-dimensional model (using two auxiliary variables to explain possible interaction of unresolved scales with these two modes). Nevertheless, the forecasting skill of these nonparametric models are relatively comparable to those of the parametric model in \cite{cmg:14}, although their testing and verification periods are slightly different. In their paper, they analyzed the skill of their model for different years. While the nonparametric model we proposed here does not give any mean to understand the physics, it is still advantageous to have the nonparametric model as a reference for developing a more physics-based parametric model. Particularly in this example, it can be concluded that the parametric model in \cite{cmg:14} is an appropriate model for predicting these two-dimensional MJO modes, since their predictive skill is relatively comparable (or not worse) compared to that of the diffusion forecasting model.

\section{Concluding discussion}\label{sec6}

In this paper, we introduced new basis functions for practical implementation of the diffusion forecasting, a nonparametric model that approximates the Fokker-Planck equation of stochastic dynamical systems. The proposed approach can avoid the expensive computation of a large number of leading eigenvectors of a diffusion matrix $T$ of size $N\times N$, where $N$ denotes the size of the training data set and could be very large. The new basis functions are constructed using the truncated unpivoting Householder QR-decomposition of an $N\times (M_E+M_Q)$ matrix, which consists of consecutive $M_Q$ columns of $T$ and it's $M_E$ leading eigenvectors. Numerically, the $QR$-decomposition of a $N\times (M_E+M_Q)$ matrix is more efficient compared to finding leading $M_E+M_Q$ eigenvectors of an $N\times N$ matrix. 

As long as $M_Q$ is large enough such that the time series corresponding to the selected columns fill up the invariant measure of the dynamical system, the new basis functions perform well in diffusion forecasting. Visually, the supports of the new basis functions gradually grow from a local regime to the whole domain; the basis function oscillates within its supports behaving like a windowed Fourier mode. The locality of the leading modes by QR factorization is inherent because of the diagonally banded structure of $T$. Using the new basis functions for the $M$-term approximation of the densities in diffusion forecasting is essentially similar to a multiresolution approximation to the densities.

The effectiveness of the proposed algorithm has been demonstrated by numerical examples in both deterministically chaotic and stochastic dynamical systems; in the former case, the superiority of the proposed basis functions over purely eigenvectors is significant, while in the latter case forecasting accuracy is improved relative to using a purely small number of eigenvectors. {\color{black} For the deterministic systems, our numerical results suggest that if the total number of basis functions used is large enough, the accurate estimation with the mixed basis can also be achieved with the purely QR-decomposed basis. Beyond the computational advantages of avoiding solving a large eigenvalue problem, we suspect that the forecast improvement over using the purely eigenvectors (even with the same number of total eigenmodes as shown in Lorenz-63 example) can be due to: 1) The densities are not necessarily smooth functions so we expect that local basis functions such as those from QR-decomposition could be better to represent non-smooth functions if they can significantly reduce the approximation error in \eqref{fouriersum}-\eqref{DF:c} while keeping the approximation error in \eqref{DF:approxA1} reasonably small. 2) As $N$ becomes large (still finite), accurate estimation of a large number of leading eigenvectors of matrix $T$ is computationally challenging. Essentially we have no access to the accurately estimated basis functions. For stochastic dynamical systems with elliptic generator as considered in this paper, the densities are smooth and therefore it is more reasonable to represent them with smooth eigenbasis which agrees with the theoretical analysis in \cite{bgh:15} and \cite{harlim:18}. This is a possible reason why we don't see improvement over using large number of eigenbasis functions. So for stochastic system, the advantage of the new basis functions is its accessibility. These findings suggest that these new basis functions can be useful for large-scale problems when only a few leading eigenfunctions can be estimated accurately.}

%
%

{\bf Acknowledgments.} The research of J.H. is partially supported by the Office of Naval Research Grants N00014-16-1-2888 and the National Science Foundation grant DMS-1317919, DMS-1619661. We thank D.Giannakis for providing the NLSA modes for the experiments in Section~\ref{sec5}. H.Y. is supported by the startup package of the Department of Mathematics, National University of Singapore.

\bibliographystyle{apsrev4-1}
\bibliography{ref}

\end{document}